\documentclass[a4paper, 12pt]{article}
\usepackage{graphicx}          
\usepackage{amsmath}
\usepackage{amssymb}                              
\def\dom{{\rm dom}}

\newcommand{\halmos}{\rule{1ex}{1.4ex}}
\newcommand{\qed}{\hfill \mbox{$\halmos$}}

\newcommand{\kk}{{\mathcal K}}
\newcommand{\ki}{{\mathcal K}_\infty}
\def \R{{\mathbb R}}
\def \N{{\mathbb N}}

\def\taus{\tau_\sigma}

\def\So{\mathcal{S}}

\def\Sd{\mathcal{S}_{d}}
\def\Sa{\mathcal{S}_{a}}

\def\T{\mathcal{T}}
\def\Ta{\mathcal{T}_a}

\def\Td{\mathcal{T}_d}
\def\U{\mathcal{U}}
\def\Ou{\mathcal{O}}

\def\V{\mathcal{V}}
\def\X{\mathcal{X}}
\def\Gr{{\rm Graph}}
\newtheorem{remark}{\itshape Remark}[section]
\newtheorem{prop}{\itshape Proposition}[section]
\newtheorem{corollary}{\itshape Corollary}[section]
\newtheorem{theorem}{\itshape Theorem}[section]
\newtheorem{lem}{\itshape Lemma}[section]
\newtheorem{as}{Assumption}
\newtheorem{definition}{\itshape Definition}[section]
\title{\LARGE \bf
Invariance Principles for Switched Systems with Restrictions}

\author{J. L. Mancilla-Aguilar \thanks{ \jma}
 \;and R.A. Garc\'{\i}a \thanks{\rg}}

\def\jma{Department of Physics \&
        Mathematics, Instituto Tecnol\'ogico de Buenos Aires
        {\tt\small jmancill@itba.edu.ar}.}

\def\rg{Department of Physics \&
        Mathematics, Instituto Tecnol\'ogico de Buenos Aires
        {\tt\small ragarcia@itba.edu.ar}.}

\begin{document}
\maketitle

\begin{abstract}                         
In this paper we consider switched nonlinear systems under average dwell time switching signals,  with
an otherwise arbitrary compact index set and with additional constraints in the  switchings. We present
invariance principles for these systems and derive by using observability-like
notions some convergence and asymptotic stability criteria. These results enable us to
analyze the stability of solutions of switched systems with both state-dependent constrained switching and switching
whose logic has memory, i.e., the active subsystem only can switch to a prescribed
subset of subsystems.
\end{abstract}

\section{Introduction}
A switched system is a family of continuous-time dynamical subsystems and a rule,
usually  time or state-dependent, that orchestrates the switching between them. At
first glance switched systems may look simple; nevertheless  their behavior may be
very complicated, being a classical example of this fact, divergent trajectories
obtained by switching among asymptotically stable subsystems  (see
\cite{Liberzon-book}). Consequently, the stability analysis of such systems turned
out to be an important and challenging problem which has received considerable
attention in the recent literature (see \cite{DeCarlo}, \cite{Liberzon-book},
\cite{Lin} and references therein). Although the stability of switched systems under
arbitrary switching laws can be assured by the existence of a common Lyapunov
function (CLF) for all the switching modes (\cite{Liberzon-book},
\cite{Mancilla-Garcia-scl00}), this type of stability condition is deemed to be too
conservative when a particular type of switching logic is considered. In fact,
switched systems that do not share a CLF may be stable under restricted switching
signals. Restrictions on the set of admissible switching signals of a certain
switched system arise naturally from physical constraints of the system, from design
strategies (e.g. discontinuous control feedback laws), or from the knowledge about
possible switching logic of the switched system, e.g., partitions of the state space
and their induced switching rules. Multiple Lyapunov functions (MLF) have been shown
to be very useful tools for the stability analysis of switched systems with
constrained switchings. In this context each switching mode may have its own
Lyapunov function (see \cite{DeCarlo} and references therein). However, some
additional conditions are necessary to assure that the value of each Lyapunov
function on its corresponding mode will decrease. Sufficient conditions for
asymptotic stability of switched systems with MLF can be found in \cite{DeCarlo},
\cite{Liberzon-book} and references therein. When the derivative of a candidate
Lyapunov function with respect to (w.r.t) each mode is only {\em non-positive}, the convergence of the solutions of the switching system to an equilibrium point, and consequently the asymptotic stability, can be  derived from one of the
various extensions to switched systems of LaSalle's invariance principle for differential equations (see \cite{LaSalle, LaSallebook}). Hespanha
in \cite{Hespanha} introduced an invariance principle for switched linear systems
under {\em persistently dwell-time} switching signals and in \cite{Hespanha2}
Hespanha {\em et al.} extended some of those results to a family of nonlinear
systems. Bacciotti and Mazzi
 presented  in \cite{Bacciotti} an invariance principle for switched
systems with {\em dwell-time} signals. An invariance principle for switched
nonlinear systems with {\em average dwell-time} signals that satisfy {\em
state-dependent} constraints was derived by Mancilla-Aguilar and Garc\'{\i}a in
\cite{Mancilla-Garcia-scl06} from the sequential compactness of particular classes
of trajectories of  switched systems. Based on invariance results for hybrid systems
(\cite{Sanfelice-Goebel-Teel}), Goebel {\em et al.} in \cite{Goebel} obtained
recently invariance results for switched systems under various types of switching
signals. Lee and Jiang in \cite{Lee} gave a generalized version of
Krasovskii-LaSalle Theorem for time-varying switched systems. Under certain
ergodicity conditions on the switching signal, some stability results were also obtained
in \cite{Cheng,Wang, Wang2}.

Most of the invariance results for switched systems already published only consider
restrictions originated by the timing of the switchings or by the state dependence of
it. Nevertheless there is also an important restriction to take into account: the
fact that not all the subsystems may be accessible from a particular one, i.e. the
case in which the switching logic has memory. This restriction is clearly exhibited,
for example, in switched systems which are the continuous portion of a hybrid
automaton (see \cite{DeCarlo}, \cite{lygeros}). In this regard, the invariance
principles developed for hybrid systems in \cite{lygeros} and
in \cite{Sanfelice-Goebel-Teel} could be useful in the analysis of switched systems
with this class of restriction in the switchings.

In this paper we present invariance results that hold for trajectories of switched
systems with a non necessarily finite number of subsystems and whose switching
signals verify an average dwell time condition and belong to a family for which a certain property {\bf P} holds. As various of
the restricted switching classes mentioned above satisfy {\bf P}, these results
enable us to obtain in an unified way invariance theorems for all of them.  Based on
these invariance results, we derive new convergence and stability criteria that
recover, generalize and strengthen some results previously obtained. In particular:
\begin{itemize}
\item Theorem \ref{thVG} extends LaSalle's invariance principle
 to switched systems with different restrictions on the
 switching signals by involving both backward and forward invariance as in, for
 example, \cite[Theorem 6.4]{LaSallebook}.
\item Theorems 1 and 2 in \cite{Bacciotti}, Corollary 5.6 of \cite{Goebel} and
Proposition 4.1 of \cite{Mancilla-Garcia-scl06} follow from Theorem \ref{theoremVG}.
\item The first conclusion of Corollary 4.10 in \cite{Goebel} is a particular case
of Theorem \ref{convergence0}. \item Theorem \ref{guas2bis} is an improvement of
Theorem 15 in \cite{Wang}. In fact the hypotheses of Theorem \ref{guas2bis} are
weaker since the existence of a Common Joint Lyapunov Function is not assumed. \item
Corollary \ref{final}, whose hypotheses are weaker than those of  Theorem 3 in
\cite{Cheng}, improves it.
\end{itemize}

The paper presents two groups of statements. First we present statements about
invariance of sets to which bounded trajectories of the switched systems converge
(Theorems \ref{thVG}, \ref{theoremVG} and \ref{thHG}). These statements involve either
continuous functions which are nonincreasing along forward complete trajectories of
the switched system or appropriately fast vanishing ``outputs''. Finally we present
results about convergence and asymptotic stability (Theorems \ref{convergence0} to
\ref{guas2bis}) that rely on observability-like conditions on the functions which
bound the derivatives of nonincreasing functions as those mentioned above.

The article unfolds as follows. Section 2. contains the basic definitions. In
Section 3. we present invariance principles for switched systems with constrained
switching. Convergence and stability results are given in Section 4. Finally
Section 5. contains some conclusions.
\section{Basic definitions}\label{basdef}
In this work we consider switched systems described by
\begin{eqnarray}
\dot{x}= f(x,\sigma) \label{ss}
\end{eqnarray}
where $x$ takes values in $\R^n$, $\sigma: \R \rightarrow \Gamma$, with $\Gamma$ a
compact metric space, is a {\em switching signal}, {\em i.e.}, $\sigma$ is piecewise
constant (it has at most a finite number of jumps in each compact interval) and is
continuous from the right and $f: \dom(f) \to \R^n$, with $\dom(f)$  a closed
subset of $\R^n \times \Gamma$, is continuous.

For each $\gamma \in \Gamma$, let $\chi_\gamma=\{\xi\in \R^n:\;(\xi,\gamma)\in
\dom(f)\}$ and $f_\gamma:\chi_\gamma\to \R^n$ be defined by
$f_\gamma(\xi)=f(\xi,\gamma)$; then $\chi_\gamma$ is closed and $f_\gamma$ is
continuous for any $\gamma \in \Gamma$. We note that when $\Gamma$ is finite, these
last two conditions imply that $\dom(f)$ is closed and that $f$ is continuous on
$\dom(f)$.  In the sequel we denote with $\So$ the set of all the switching signals.

Given $\sigma \in \So$, a solution of (\ref{ss}) corresponding to $\sigma$ is a
locally absolutely continuous function $x:I_x \to \R^n$, with $I_x \subset \R$ a
nonempty interval, such that $(x(t),\sigma(t))\in \dom(f)$ for all $t \in I_x$ and
$\dot{x}(t)= f(x(t),\sigma(t))$ for almost all $t \in I_x$. The solution $x$ is
complete if $I_x=\R$ and forward complete if $\R_{\ge 0} \subset I_x$. A pair
$(x,\sigma)$ is a trajectory of (\ref{ss}) if $\sigma \in \So$ and $x$ is a solution
of (\ref{ss}) corresponding to $\sigma$. The trajectory is complete or forward
complete if $x$ is complete or forward complete, respectively. Given a subset $\Ou$
of $\R^n$, we say that the trajectory $(x,\sigma)$ is precompact relative to $\Ou$
if there exists a compact set $B \subset \Ou$ such that $x(t)\in B$ for all $t \in
I_x$. When $\Ou=\R^n$ we simply say that $(x,\sigma)$ is precompact.
\begin{remark}
Note that we do not suppose that $\dom(f)=\R^n \times \Gamma$. In this way
we can take into account, in the analysis of the asymptotic behavior of a given
trajectory $(x,\sigma)$ of (\ref{ss}), some kind of state-dependent constraints
which the trajectory under study must satisfy. In fact, in some situations we are
not interested in the behavior of an arbitrary forward complete trajectory $(x,\sigma)$ of
a switched system (\ref{ss}) (with $\dom(f)=\R^n \times \Gamma$) but only of one of those that verify the constraint
\begin{eqnarray} \label{sdc}
x(t) \in \chi_{\sigma(t)}\quad \makebox{for all}\; t \in I_x,\end{eqnarray} where
$\{\chi_\gamma:\gamma \in \Gamma\}$ is a collection of subsets of $\R^n$. If we
consider the map $\tilde{f}$, which is the restriction of $f$ to the set
$\dom(\tilde{f})=\{(\xi,\gamma):\;\xi \in \chi_\gamma\}$, and if $\dom(\tilde{f})$
is closed in $\R^n \times \Gamma$, then the set of trajectories $(x,\sigma)$ of
(\ref{ss}) which verify (\ref{sdc}) coincides with the set of trajectories of
\begin{eqnarray} \label{ssc}
\dot{x}= \tilde{f}(x,\sigma).\end{eqnarray} It must be pointed out that in this way
we can consider the system as if its switching is state-independent, and focus on the
restrictions imposed to it by the timing of the discontinuities of the switching
signal and/or by the accessibility  to certain subsystems from another ones.
\end{remark}

In this paper we consider forward complete solutions of (\ref{ss}) corresponding to
switching signals $\sigma$ which belong to particular subclasses of $\So$. Let
$\Lambda(\sigma)$ be the set of times where $\sigma$ has a jump (switching time).
Following \cite{Hespanha} we say that $\sigma \in \So$ has a dwell-time $\tau_D>0$
if $|t-t^\prime| \ge \tau_D$ for any pair $t,t^\prime \in \Lambda(\sigma)$ such that
$t\neq t^\prime$.

A switching signal $\sigma$ has an average dwell-time $\tau_D>0$ and a chatter bound
$N_0 \in \N$ if the number of switching times of $\sigma$ in any open finite
interval $(\tau_1,\tau_2) \subset \R$ is bounded by $N_0+(\tau_2-\tau_1)/\tau_D$,
i.e. ${\rm card}(\Lambda(\sigma)\cap (\tau_1,\tau_2)) \le
N_0+(\tau_2-\tau_1)/\tau_D$.

We denote by $\Sa[\tau_D,N_0]$ the set of all the switching signals which have an
 average dwell-time $\tau_D>0$ and a chatter bound $N_0 \in \N$ and by
$\Ta[\tau_D,N_0]$ the set of all the complete trajectories $(x,\sigma)$ of
(\ref{ss}) with $\sigma \in \Sa[\tau_D,N_0]$ and let $\Sa=\bigcup_{\tau_D>0, N_0> 0}
\Sa[\tau_D,N_0]$ and $\Ta=\bigcup_{\tau_D>0, N_0> 0} \Ta[\tau_D,N_0]$. We note that
the set of switching signals $\sigma$ which have a dwell-time $\tau_D>0$ coincides
with $\Sa[\tau_D,1]:=\Sd[\tau_D]$. We denote by $\Td[\tau_D]$ the set of all the
complete trajectories $(x,\sigma)$ of (\ref{ss}) with $\sigma \in \Sd[\tau_D]$ and
let $\Sd=\bigcup_{\tau_D >0} \Sd[\tau_D]$ and $\Td=\bigcup_{\tau_D>0} \Td[\tau_D]$.

For $\Gamma$ a finite set and $T>0$, we denote by $\So_e[T]$ the family of all the
switching signals $\sigma$ which verify the following ``ergodicity'' condition: for
every $t_0\ge 0$ and every $\gamma \in \Gamma$, $\sigma^{-1}(\gamma) \cap
[t_0,t_0+T]\neq \emptyset$.

$\T_e[T]$ will denote the set of complete trajectories $(x,\sigma)$ with $\sigma \in
\So_e[T]$ and $\So_e=\bigcup_{T >0} \So_e[T]$ and $\T_e=\bigcup_{T>0} \T_e[T]$.

The families of switching signals already introduced have no restrictions on the
accessibility from any subsystem to another. The family of switching signals ---and
their corresponding trajectories--- that we introduce next, takes into account the
case in which the switching logic has memory, i.e. when a subsystem corresponding to
an index $\gamma \in \Gamma$ can only switch to subsystems corresponding to modes
$\gamma^\prime$ that belong to a certain subset $\Gamma_\gamma\subset \Gamma$.

Given a set-valued map $H:\Gamma \rightsquigarrow \Gamma$, $\So^H$ is the set of all
the switching signals $\sigma$ which verify the condition $\sigma(t)\in
H(\sigma(t^-))$ for every time $t \in \Lambda(\sigma)$. Here $\sigma(t^-)=\lim_{s\to
t^-}\sigma(s)$. $\T^H$ denotes the set of all the complete trajectories $(x,\sigma)$
with $\sigma \in \So^H$.  This class of switching signals enable us, for example, to
model the restrictions imposed by the discrete process of a hybrid system whose
continuous portion is as in (\ref{ss}) (see \cite{DeCarlo}).

\section{Invariance results for trajectories which satisfy a dwell-time condition}\label{inva}
In this section we present some invariance results that enable us to characterize
the asymptotic behavior of a precompact forward complete trajectory $(x,\sigma)$ of
(\ref{ss}) with $\sigma$ belonging to a certain subclass of $\Sa$. The consideration
of such subclass allows us to obtain in an unified way invariance results for
systems whose switching signals  undergo different restrictions.

We recall that a point $\xi \in \R^n$ belongs to $\Omega(x)$, the $\omega$-limit set
of $x:I_x \rightarrow \R^n$, with $\R_{\ge 0} \subset I_x$, if there exists a
strictly increasing sequence of times $\{s_k\} \subset I_x$ with $\lim_{k
\rightarrow \infty}s_k= + \infty$ and $\lim_{k \rightarrow \infty}x(s_k)=\xi$. The
$\omega$-limit set $\Omega(x)$ is always closed and, when $x$ evolves in a compact
set of $\R^n$, it is nonempty, compact, connected if $x$ is continuous,  and $x \to
\Omega(x)$ (for a set $M \subset \R^n$, $x \to M$ if $\lim_{t \rightarrow +\infty}
d(x(t),M)=0$, being $d(\xi,M)=\inf_{\nu \in M}|\nu-\xi|$).

As was done in \cite{Mancilla-Garcia-scl06}, we will associate to each 
forward complete trajectory $(x,\sigma)$ of (\ref{ss}) with $\sigma \in \Sa$, 
the nonempty set $\Omega^{\sharp}(x,\sigma)\subset \R^n \times \Gamma$ that we
introduce in the following
\begin{definition} \label{omega} Given a forward complete trajectory $(x,\sigma)$ of
(\ref{ss}) with $\sigma \in \Sa$, a point $(\xi, \gamma)\in \R^n \times \Gamma$
belongs to $\Omega^\sharp(x,\sigma)$ if there exists a strictly increasing and
unbounded sequence $\{s_k\} \subset \R_{\ge 0}$ such that
\begin{enumerate}
\item $\lim_{k \rightarrow \infty} \taus^1(s_k)-s_k = r, \; 0<\,r\,\leq \infty$,
\item $\lim_{k\rightarrow \infty}x(s_k)=\xi$ and $\lim_{k \rightarrow
\infty}\sigma(s_k)=\gamma$.
\end{enumerate}
Here, for any $t\in \R$, $\taus^1(t)=\inf\{s \in \Lambda(\sigma):t < s\}$ if $\{s
\in \Lambda(\sigma):t < s\}\neq \emptyset$ and $\taus^1(t)=+\infty$ in other case
(i.e. $\taus^1(t)$ is the first switching time greater than $t$).
\end{definition}

Let $\pi_1:\R^n \times \Gamma \rightarrow \R^n$ be the projection onto the first
component. Then the  following relation between $\Omega(x)$ and
$\Omega^\sharp(x,\sigma)$ holds.
\begin{lem} \label{pi1}  Let
$(x,\sigma)$ be a forward complete trajectory of (\ref{ss}) with $\sigma \in \Sa$
that is precompact relative to $\Ou \subset \R^n$.
 Then $\Omega^\sharp(x,\sigma)\subset \dom(f) \cap (\Ou \times \Gamma)$ and $\Omega(x)=\pi_1(\Omega^\sharp(x,\sigma))$.
\end{lem}
{\bf Proof.}
That $\Omega^\sharp(x,\sigma)\subset \dom(f)\cap (\Ou \times \Gamma)$ follows from
the fact that for all $t \in I_x$ $(x(t),\sigma(t))$ belongs to a compact subset of
$\dom(f) \cap (\Ou \times \Gamma)$  and from the definition of
$\Omega^\sharp(x,\sigma)$.  The proof of the other assertion follows
 {\em mutatis mutandis} from the proof of Lemma 4.1 in \cite{Mancilla-Garcia-scl06}. \qed

In order to see that the set $\Omega^\sharp(x,\sigma)$ enjoys certain kind of
invariance property, let us introduce the following
\begin{definition} Given a family $\T^*$ of complete trajectories of (\ref{ss}), we say that a
nonempty subset $M \subset \R^n \times \Gamma$ is {\em weakly-invariant w.r.t}
$\T^*$ if for each $(\xi,\gamma)\in M$ there is a trajectory $(x,\sigma) \in \T^*$
such that $x(0)=\xi$, $\sigma(0)=\gamma$ and $(x(t),\sigma(t)) \in M$ for all $t \in
\R$.
\end{definition}

This notion of weak invariance differs from the one introduced in
\cite{Mancilla-Garcia-scl06}, in that the last one involves only forward invariance
while the introduced here also involves backward invariance.

Let us introduce now the following class of switching signals.
\begin{definition} \label{P} We say that a family of switching signals $\So^*$
has the property {\bf P} if
\begin{enumerate}
\item  $\So^*\subset \Sa[\tau_D,N_0]$ for some $\tau_D>0$ and some $N_0\in \N$;
\item for any $s>0$ and any $\sigma \in \So^*$, $\sigma(\cdot+s)\in \So^*$; \item
for every sequence $\{\sigma_k\} \subset \So^*$, there exist $\sigma^*\in \So^*$ and
a subsequence $\{\sigma_{k_l}\}$ such that $\lim_{l\to
\infty}\sigma_{k_l}(t)=\sigma^*(t)$ for almost all $t\in \R$.
\end{enumerate}
\end{definition}
\begin{lem}\label{unif1} The following classes of switching signals have the
property {\bf P}:
\begin{enumerate}
\item $\Sa[\tau_D,N_0]$ for every $\tau_D>0$ and every $N_0 \in \N$; \item
$\Sd[\tau_D]\cap \So^H$ for all $\tau_D>0$ and every $H :\Gamma \rightsquigarrow
\Gamma$ such that the set $\Gr(H)=\{(\gamma,\gamma^\prime)\in \Gamma \times
\Gamma:\;\gamma^\prime \in H(\gamma)\}$ is closed; \item $\Sd[\tau_D]\cap \So_e[T]$
for every $\tau_D>0$ and every $T>0$.
\end{enumerate}
\end{lem}
{\bf Proof.} See Appendix. \qed

The next result will be instrumental in what follows.
\begin{theorem}\label{ginvariance} Let $\So^*$ be a family of switching
signals which verifies property {\bf P} and let $\T^*$ be the set of all the
complete trajectories $(\overline{x},\overline{\sigma})$ of (\ref{ss}) with $\overline{\sigma} \in \So^*$. Then, if
$(x,\sigma)$ is a precompact forward complete trajectory of (\ref{ss}) such that
$\sigma \in \So^*$, $\Omega^\sharp(x,\sigma)$ is weakly-invariant w.r.t $\T^*$.
\end{theorem}
{\bf Proof.} See Appendix. \qed
\begin{remark}\label{unif2}
Since the weak invariance of $\Omega^\sharp(x,\sigma)$ is a cornerstone of the
results that we present below (Theorems \ref{thVG} to \ref{thHG}), Theorem
\ref{ginvariance} enables us to obtain in a unified way invariance results not only
for the different switching signals explicitly mentioned in Lemma \ref{unif1}  but
also for any other that verifies property {\bf P}.
\end{remark}
\begin{remark}
 At first glance, it would seem more natural to consider for a given precompact forward complete
 trajectory $(x, \sigma)$ of (\ref{ss}) its $\omega$-limit set $\Omega(x, \sigma)$ instead
 of $\Omega^\sharp(x,\sigma) \subset
 \Omega(x,\sigma)$. Nevertheless, there exist forward complete trajectories $(x, \sigma)$ of (\ref{ss}) with $\sigma \in
\Sa$ such that $\Omega(x,\sigma)$ is not weakly-invariant for any family of
trajectories of that switched system.
\end{remark}

Next, we present two invariance results that involve the existence of a function $V$
which is nonincreasing along a trajectory of (\ref{ss}). In order to do so, we
introduce the following class of functions.

\begin{definition}\label{claseV}
We say that a function $V : \dom(V) \rightarrow \R$ belongs to class $\mathcal{V}$,
if it verifies
\begin{enumerate}
 \item $\dom(V) \subset \R^n \times \Gamma$.
\item For every $\gamma \in \Gamma$, $\mathcal{D}_\gamma :=\{\xi \in \R^n : (\xi,
\gamma) \in \dom(V)\}$ is an open set. \item Let $\Ou := \pi_1(\dom(V))$. Then
$\Ou_\gamma := \Ou \cap \chi_\gamma \subset \mathcal{D}_\gamma \; \forall \gamma \in
\Gamma$ \item For all $\gamma \in \Gamma$, $V_\gamma(\cdot) := V(\cdot, \gamma)$ is
differentiable on $\Ou_\gamma$.
\end{enumerate}
\end{definition}
We note that $\dom(f)\cap (\Ou \times \Gamma)=\cup_{\gamma \in
\Gamma}(\Ou_\gamma \times \{\gamma\})\subset \dom(f)\cap \dom(V)$.

We also note that when $\Gamma$ is finite, the restriction of any function $V \in
\mathcal{V}$ to $\dom(f) \cap (\Ou \times \Gamma)$ is continuous.

 In what follows, for a function $V \in \mathcal{V}$, let
$Z_V=\{(\xi,\gamma)\in \dom(f)\cap (\Ou \times \Gamma):\;\nabla V_\gamma(\xi)
f_\gamma(\xi)=0\}$. 
\begin{as} \label{V} The forward complete trajectory $(x,\sigma)$ of
(\ref{ss}) verifies the following: there exists a function $V \in \mathcal{V}$ whose
restriction to $\dom(f) \cap (\Ou \times \Gamma)$ is continuous,  $(x, \sigma)$ is
precompact relative to $\Ou$ and $v(t)=V(x(t),\sigma(t))$ is nonincreasing on $[0, +
\infty)$.
\end{as}

\begin{theorem}\label{thVG} Let $\So^*$ be a family of switching signals which has property {\bf P}
and let $\T^*$ be the set of all the complete trajectories $(x,\sigma)$ of
(\ref{ss}) with $\sigma \in \So^*$. Suppose that $(x,\sigma)$, with $\sigma \in
\So^*$,  is a forward complete trajectory of (\ref{ss}) for which Assumption \ref{V}
holds. Then there exists $c\in \R$ such that $x \to \pi_1(M(c))$, where $M(c)$ is
the maximal weakly-invariant set w.r.t. $\T^*$ contained in $V^{-1}(c) \cap Z_V$.
\end{theorem}

{\bf Proof.} Since $\Omega^\sharp(x,\sigma)$ is weakly-invariant
w.r.t. $\T^*$ and, from Lemma \ref{pi1}, $x\to \pi_1(\Omega^\sharp(x,\sigma))$, we
only have to prove that $\Omega^{\sharp}(x,\sigma)\subset V^{-1}(c)\cap Z_V$ for
some $c\in \R$.

As $(x,\sigma)$ is precompact relative to $\Ou$, there exists a compact set $B
\subset \Ou$ such that $x(t) \in B$ for all $t \in I_x$. Therefore
$(x(t),\sigma(t))$ belongs to the compact set $\dom(f)\cap (B \times \Gamma)$ for
all $t \in I_x$. Thus $v(t)$ is bounded, since $V$ is continuous on $\dom(f) \cap (B
\times \Gamma)$, and nonincreasing by hypothesis; in consequence there exists
$\lim_{t \to +\infty}v(t)=c$.

Let $(\xi,\gamma) \in \Omega^{\sharp}(x,\sigma)$. Then there exists a 
strictly increasing and unbounded sequence $\{s_k\}$ which verifies 1. and 2. of
Definition \ref{omega}. Since $(x(s_k),\sigma(s_k)) \to (\xi,\gamma)$ as $k  \to
\infty$,  $c=\lim_{k\to \infty}v(s_k)=\lim_{k \to \infty}V(x(s_k),\sigma(s_k))$ $ =
V(\xi,\gamma)$ and $(\xi,\gamma) \in V^{-1}(c)$. Let us show that $(\xi,\gamma)$
also belongs to $Z_V$.

As $\Omega^{\sharp}(x,\sigma)$ is weakly-invariant w.r.t. $\T^*$ there exists
$(x^*,\sigma^*)\in \T^*$ such that $(x^*(0),\sigma^*(0))=(\xi,\gamma)$ and
$(x^*(t),\sigma^*(t)) \in \Omega^{\sharp}(x,\sigma)$ for all $t \in \R$. Then,
taking into account that $\Omega^{\sharp}(x,\sigma)\subset V^{-1}(c)$,
$V(x^*(t),\sigma^*(t))=c$ for all $t \in \R$. In particular, since
$\sigma^*(t)=\gamma$ on $[0,\tau)$ for $\tau$ small enough, then $V(x^*(t),
\gamma)=c$ on $[0,\tau)$. Therefore $\nabla V_\gamma(\xi)f_\gamma(\xi)=0$, and
$(\xi,\gamma) \in Z_V$.\qed
\begin{remark} We note that Theorem \ref{thVG} is an extension to
switched systems of the well known LaSalle's invariance principle for 
differential equations (see, for example, \cite[Theorem 6.4]{LaSallebook}).
\end{remark}

In the sequel, for any $\sigma \in \So$ and any $\gamma \in \Gamma$,  let
$\mathcal{I}_{\sigma, \gamma} = \sigma^{-1}(\gamma) \cap [0, + \infty)$.

When  $\Gamma$ is a finite set, we can relax the nonincreasing condition in
Assumption \ref{V} as follows.

\begin{as} \label{MV} The forward complete trajectory $(x,\sigma)$ of
(\ref{ss}) verifies the following:  there exists a function $V \in \mathcal{V}$ such
that  $(x, \sigma)$ is precompact relative to $\Ou$ and $v(t)=V(x(t),\sigma(t))$ is
nonincreasing on $\mathcal{I}_{\sigma, \gamma}$, for all $\gamma \in \Gamma$.
\end{as}
\begin{remark}\label{fini}\mbox{}
 Assumptions of this kind are standard when the stability analysis of
the zero solution of a switched system is performed by means of multiple Lyapunov
functions (see \cite{DeCarlo}, \cite{Liberzon-book}).
\end{remark}

In what follows, when $\Gamma$ is a finite set, we identify it with the set
$\{1,\ldots,N\}$ $\subset \N$, where $N={\rm card}(\Gamma)$.
\begin{theorem}\label{theoremVG}  Suppose that $\Gamma$ is finite and let $\So^*$ and $\T^*$ be as in
Theorem \ref{thVG}. Suppose that $(x,\sigma)$, with $\sigma \in \So^*$, is
 a forward complete trajectory of (\ref{ss}) for which Assumption
\ref{MV} holds.  Then there exists $\vec{c} = (c_1, \ldots, c_N) \in \R^N$ such that $x \to
\pi_1(M(\vec{c}))$, where $M(\vec{c})$ is the maximal weakly-invariant set w.r.t.
$\T^*$ contained in $\cup_{\gamma \in \Gamma} \{(\xi,\gamma)\in \dom(f) \cap (\Ou
\times \Gamma):\; V_\gamma(\xi)=c_\gamma \} \cap Z_V$.
\end{theorem}

{\bf Proof.} For $\gamma \in \Gamma$ we define $c_\gamma$ as
follows:
\begin{enumerate}
\item $c_\gamma=\lim_{t\to +\infty, \: t \in \mathcal{I}_{\sigma, \gamma}} v(t)$ if
$\mathcal{I}_{\sigma, \gamma}$ is unbounded. (This limit exists since $v$ is
non-increasing and bounded on $\mathcal{I}_{\sigma, \gamma}$). \item $c_\gamma=a$
for every $\gamma$ such that $\mathcal{I}_{\sigma, \gamma}$ is bounded. Here $a \in
\R$ is arbitrary.
\end{enumerate}
Reasoning as in the proof of Theorem \ref{thVG}, in order to prove the thesis it
suffices to show that $\Omega^{\sharp}(x,\sigma) \subset \cup_{\gamma \in \Gamma}
\{(\xi,\gamma)\in \dom(f)\cap (\Ou \times \Gamma):\;V_\gamma(\xi)=c_\gamma\}\cap
Z_V$.

Let $(\xi,\gamma) \in \Omega^{\sharp}(x,\sigma)$. Then there exists a strictly increasing and
unbounded sequence $\{s_k\}$ which verifies 1. and 2. of Definition \ref{omega}.
Since $\sigma(s_k)\to \gamma$ and $\Gamma$ is a finite set, $\sigma(s_k)=\gamma$ for
$k$ large enough and for those $k$, $s_k \in \mathcal{I}_{\sigma,\gamma}$. It
follows that $V(x(s_k),\gamma) \to c_{\gamma}$ as $k \to \infty$ and in consequence,
$V(\xi,\gamma)=c_{\gamma}$ and $(\xi,\gamma) \in \cup_{\gamma^\prime \in \Gamma}
\{(\xi^\prime,\gamma^\prime)\in \dom(f)\cap (\Ou \times
\Gamma):\;V(\xi^\prime,\gamma^\prime)=c_{\gamma^\prime}\}$. That
$(\xi,\gamma)\in Z_V$ can be proved in the same way as
 in the proof of Theorem \ref{thVG}.\qed
\begin{remark} Some invariance results for switched systems reported in the literature can
be derived from Theorem \ref{theoremVG}. In particular \cite[Theorems 1 and
2]{Bacciotti}, \cite[Proposition 4.1]{Mancilla-Garcia-scl06} and \cite[Corollary
5.6]{Goebel}.
\end{remark}

 The following invariance result involves {\em weakly meagre functions}. We
recall that a function $y:\R_{\ge 0}\to \R$ is weakly meagre if $\lim_{k \to
\infty}(\inf\{|y(t)|:t\in I_k\})=0$ for every family $\{I_k:\;k \in \N\}$ of
nonempty and pairwise disjoint intervals in $\R_{\ge 0}$ with $\inf\{\mu(I_k):\;k\in
\N\}>0$, where $\mu$ stands for the Lebesgue measure (see \cite{Logemann-Ryan}). We
note that, for example, any function $y \in L^{p}([0,\infty))$ with $p > 0$ is
weakly meagre. More generally, if there exist positive numbers $\tau$ and $p$  such
that $\int_{t}^{t+\tau} |y(s)|^p ds$ converges to $0$ as $t \to +\infty$, then $y$
is weakly meagre.
\begin{theorem} \label{thHG} Let $\So^*$ and $\T^*$ be as in
Theorem \ref{thVG}. Suppose that $(x,\sigma)$, with $\sigma \in \So^*$, is
 a forward complete trajectory of (\ref{ss}). Suppose in addition that there exists a continuous
function $h:\dom(h) \to \R$,  with $\dom(h) \subset \R^n \times \Gamma$, such that $(x(t), \sigma(t))$ evolves in a compact
subset $K$ of $\dom(h)$ for all $t \geq 0$ and that $y(\cdot)=h(x(\cdot),\sigma(\cdot))$
is weakly meagre. Then $x \to \pi_1(M^*)$, where $M^*$ is the maximal weakly
invariant set w.r.t. $\T^*$ contained in $h^{-1}(0) \cap {\rm dom}(f)$.
\end{theorem}

{\bf Proof.} Since $(x(t), \sigma(t))$ evolves in the compact
subset $K$ of $\dom(h)$ for all $t\ge 0$, we have that $(x,\sigma)$ is precompact
and that $\Omega^\sharp(x,\sigma)\subset \dom(h) \cap \dom(f)$. By similar considerations as those in the proofs of
the previous invariance results, it suffices to show that
$\Omega^\sharp(x,\sigma)\subset h^{-1}(0)$.

Let $(\xi^*,\gamma^*) \in \Omega^{\sharp}(x,\sigma)$. Then there exists a strictly
increasing and unbounded sequence $\{s_k\}$
 which verifies 1. and 2. of 
Definition \ref{omega} with $(\xi^*,\gamma^*)$ instead of $(\xi,\gamma)$. We can
assume that $\tau^1_{\sigma}(s_k)-s_k \ge 3r/4$ for all $k$.

We will construct by recursion a sequence of times $\{s^*_m\}$ and a subsequence
$\{s_{k_m}\}$ of $\{s_k\}$ such that $s_{k_m} \le s^*_m\le s_{k_m}+2^{-m}r$  and
$|y(s^*_m)|\le 1/m$ for all $m$.

For any $m \in \N$, let $r_m=2^{-m}r$ and $I_k^m=[s_k,s_k+r_m]$ for all $k \in \N$.
Since $h$ is weakly meagre, $\lim_{k \to \infty}(\inf\{|y(t)|:t\in I^1_k\})=0$.
Thus, there exist $k^*\in \N$ an a time $t \in I^1_{k^*}$ such that $|y(t)|\le 1$.
Let $s^*_1=t$ and $k_1=k^*$.

Suppose that we have already defined $\{s^*_m\}_{m=1}^l$ and $\{s_{k_m}\}_{m=1}^l$.
As $h$ is weakly meagre then $\lim_{k \to \infty}(\inf\{|y(t)|:t\in I^{l+1}_k\})=0$,
and hence there exist $k^\prime \in \N$, with $k^\prime > k_l$, and $t^\prime \in
I^{l+1}_{k^\prime}$ such that $|y(t^\prime)|\le 1/(l+1)$. We define
$s^*_{l+1}=t^\prime$ and $k_{l+1}=k^\prime$.

Consider now the sequence $\{(x(s^*_m),\sigma(s^*_m))\}$; since $s^*_m \in
[s_{k_m},s_{k_m}+r_m]$, $r_m \le r/2$ and $ \tau^1_{\sigma}(s_{k_m})-s_{k_m}\ge
3r/4$, then for all $m,\, \sigma(s^*_m)= \sigma(s_{k_m})$  and hence
$\sigma(s^*_m)\to \gamma^*$.

Given that for every $t \ge 0,\, (x(t),\sigma(t))$ belongs to the compact set
$\dom(f)\cap K$ and since $f$ is continuous on that set, there exists $M\ge 0$ such
that $|\dot{x}(t)|\le M$ a.e. on $I_x$ and in consequence, $|x(t)-x(s)|\le M|t-s|$
for all $t,s\in I_x$. Therefore $|x(s^*_m)-x(s_{k_m})|\le r M 2^{-m}$ for all $m$.
Taking into account that $x(s_{k_m})\to \xi^*$, we have that $x(s^*_{m})\to \xi^*$.

Finally, since by construction $h(x(s^*_{m}),\sigma(s^*_m))=y(s^*_{m})\to 0$ and $h$
is continuous, then $h(\xi^*,\gamma^*)=0$. \qed

\section{Convergence and stability results}\label{conve}
In this section we derive, from the invariance principles presented in \S
\ref{inva}, some convergence and stability results for switched systems with
constrained switchings.
\subsection{Convergence results}
 Let us first introduce some observability-like definitions.

Given a subset $\X \subset \R^n$, a continuous map $g:\X \to \R^n$ and
 a function $h:\X\to \R$,  we say that
for a given $\tau$ ($\tau>0$ or $\tau=\infty$) a point $\xi \in \X$ belongs to the
set $\X^f(g,h,\tau)$ (resp. $\X^b(g,h,\tau)$) if there exists a solution
$\varphi:[0,\tau]\to \X$ (resp. $\varphi:[-\tau,0] \to \X$) of $\dot{x}=g(x)$ such
that $\varphi(0)=\xi$ and $h(\varphi(t))=0$ for all $t \in [0,\tau]$ (resp. $t\in
[-\tau,0]$).

Let also  the sets $\X^f(g,h)=\bigcup_{\tau
> 0}\X^f(g,h,\tau)$, $\X^b(g,h)=\bigcup_{\tau > 0}\X^b(g,h,\tau)$
and $\X(g,h)=\X^f(g,h)\cup \X^b(g,h)$.

\begin{remark}\label{obs}\makebox[.1in]{}
 \begin{enumerate}
\item The set $\X^f(g,h,\infty)$ ($\X^b(g,h,\infty)$) coincides with the maximal
weakly forward(backward) invariant set w.r.t. $g$ contained in the set $\{\xi \in
\X:\:h(\xi)=0\}$.

We recall that a subset $K \subset \R^n$ is weakly forward(backward) invariant w.r.t
$g$ if for each $\xi \in K$ there exists a solution $\varphi:[0,\infty) \to \R^n$
($\varphi:(-\infty,0] \to \R^n$) of $\dot{x}=g(x)$ such that $\varphi(0)=\xi$ and
$\varphi(t) \in K$ for all $t \ge 0$ ($t \le 0$). \item If we consider the system
with outputs $\dot{x}=g(x)$, $y=h(x)$ and state space $\X$, with $0\in \X$, $g(0)=0$
and $h(0)=0$, then the set $\X^f(g,h)$ coincides with the set of states $\xi$ that
cannot be instantaneously distinguished from the zero state through the output $y$.
In the particular case in which $g$ is a linear function, i.e.,
$g(\xi)=A\xi$ and $h(\xi)=\xi^T C^T C\xi$, and $C$ is a matrix, then $\X(g,h)\subset \X
\cap \mathcal{U}$, being $\mathcal{U}$ the unobservable subspace of $(C,A)$. \item
When $g$ and $h$ are smooth functions we have that
$$\X(g,h) \subset \{\xi \in \X:\:L_g^kh(\xi)=0\: \forall k \in
\N_0\},$$ with $L_g^kh$ the $k$-th. Lie derivative of $h$ along $g$.
\end{enumerate}
\end{remark}
Let us introduce the following assumptions, in order to obtain some convergence
criteria based on the invariance results given in \S \ref{inva} and on the
observability-like notions already introduced.

\begin{as} \label{W} For the forward complete trajectory $(x,\sigma)$ of
(\ref{ss}) there exist a function $V \in \V$ and a family of functions
$\{W_\gamma:\Ou_\gamma \to \R, \gamma \in \Gamma\}$ such that $(x, \sigma)$ and $V$
satisfy Assumption \ref{V} and in addition 
\begin{eqnarray} \label{des1} -\nabla V_\gamma (\xi)f_\gamma(\xi)\ge W_\gamma(\xi) \ge 0 \quad
\forall \xi \in \Ou_\gamma, \quad \forall \gamma \in \Gamma.
\end{eqnarray}
\end{as}

\begin{as} \label{W2}
For the forward complete trajectory $(x,\sigma)$ of (\ref{ss}) there exist a
function $V \in \V$ which is bounded on compact subsets of $\dom(f)\cap(\Ou \times
\Gamma)$ and a family of functions $\{W_\gamma:\Ou_\gamma \to \R, \gamma \in
\Gamma\}$ such that $(x, \sigma)$ is precompact relative to $\Ou$,
$v(t)=V(x(t),\sigma(t))$ is nonincreasing on $[0, + \infty)$, (\ref{des1}) holds      and in addition
 $W(\xi,\gamma) := W_\gamma(\xi)$ is continuous on $\dom(W) = \dom(f) \cap (\Ou \times \Gamma)$.
\end{as}
\begin{as} \label{Wf}
For the forward complete trajectory $(x,\sigma)$ of (\ref{ss}) there exist a
function $V \in \V$ and a family of functions $\{W_\gamma:\Ou_\gamma \to \R, \gamma
\in \Gamma\}$ such that $(x, \sigma)$ and $V$ satisfy Assumption \ref{MV} and in
addition (\ref{des1}) holds.
\end{as}
\begin{theorem} \label{convergence0} Let $(x,\sigma)$ be a forward complete
trajectory of (\ref{ss}) with $\sigma \in \Sa$. Then the following holds:
\begin{enumerate}
 \item if $(x, \sigma)$ verifies Assumption \ref{W}, then there exists $c \in \R$ such that
$$x \to
\bigcup_{\gamma, \gamma^\prime \in \Gamma} \left( \Ou_\gamma^f(f_\gamma, W_\gamma)\cap \Ou_{\gamma^\prime}^b(f_{\gamma^\prime},W_{\gamma^\prime}) \cap
V^{-1}_\gamma(c) \cap
V^{-1}_{\gamma^\prime}(c)\right);$$ \item if $(x, \sigma)$ verifies Assumption \ref{W2}, then
$$x \to
\bigcup_{\gamma, \gamma^\prime \in \Gamma} \left( \Ou_\gamma^f(f_\gamma, W_\gamma)\cap \Ou_{\gamma^\prime}^b(f_{\gamma^\prime},W_{\gamma^\prime})\right);$$
 \item if $\Gamma$ is finite and $(x, \sigma)$ verifies Assumption \ref{Wf}, then there exists $\vec{c} \in \R^N$ such that
$$x \to
\bigcup_{\gamma, \gamma^\prime \in \Gamma} \left( \Ou_\gamma^f(f_\gamma, W_\gamma)\cap \Ou_{\gamma^\prime}^b(f_{\gamma^\prime},W_{\gamma^\prime}) \cap
V^{-1}_\gamma(c_\gamma) \cap
V^{-1}_{\gamma^\prime}(c_{\gamma^\prime})\right);$$ 

\end{enumerate}
\end{theorem}
{\bf Proof.}  Since $\sigma \in \Sa$, there exist $\tau_D>0$ and
$N_0 \in \N$ such that $\sigma \in \Sa[\tau_D,N_0]$.

In order to prove 1, let $(x,\sigma)$ verify Assumption \ref{W}. As  $\Sa[\tau_D,N_0]$ has property {\bf P} and $(x,\sigma)$ verifies the
hypotheses of Theorem \ref{thVG}, there exists a real number $c$ such that $x \to
\pi_1(M(c))$, where $M(c)$ is the maximal weakly-invariant set w.r.t.
$\Ta[\tau_D,N_0]$ contained in $V^{-1}(c)\cap Z_V$.

Let $\xi \in \pi_1(M(c))$ and $\gamma \in \Gamma$ such that $(\xi,\gamma)\in M(c)$. From the weak invariance of $M(c)$ w.r.t. $\Ta[\tau_D,N_0]$, there
 exists a trajectory $(x^*,\sigma^*) \in \Ta[\tau_D,N_0]$ such that $(x^*(0),\sigma^*(0))=
 (\xi,\gamma)$ and such that for every $t \in \R, \;\nabla V_{\sigma^*(t)}(x^*(t)) f_{\sigma^*(t)}(x^*(t)) = 0$ and $V_{\sigma^*(t)}(x^*(t))= c$. 
Since for every $\gamma \in \Gamma$, $V_\gamma(\cdot)$ is continuous on $\Ou_\gamma$, we also have that  $V_{\sigma^*(t^-)}(x^*(t))= c$ and $\nabla V_{\sigma^*(t^-)}(x^*(t)) f_{\sigma^*(t^-)}(x^*(t)) = 0$ for all $t \in \R$.
In view of (\ref{des1}),  $W_{\sigma^*(t)}(x^*(t))= W_{\sigma^*(t^-)}(x^*(t)) = 0$ for all $t\in \R$. Let us consider two cases.

{\em Case 1.}  $0 \notin \Lambda(\sigma^*)$. Then,  there exist $\tau_1 < 0 < \tau_2$ such that $\sigma^*(t) = \gamma$ for every $t \in [\tau_1, \tau_2]$. Hence, $\xi = x^*(0) \in  \Ou_\gamma^f(f_\gamma, W_\gamma, \tau_2)\cap \Ou_{\gamma}^b(f_{\gamma},W_{\gamma}, -\tau_1) \cap
V^{-1}_\gamma(c)$.

{\em Case 2.} $0 \in \Lambda(\sigma^*)$.  Let $\gamma^\prime = \sigma^*(0^-)$; then, there exist $\tau_1 < 0 < \tau_2$ such that $\sigma^*(t) = \gamma^\prime$ for all $t \in [\tau_1,0)$  and $\sigma^*(t) = \gamma$ for all $t \in [0, \tau_2]$. In consequence, $\xi = x^*(0) \in  \Ou_\gamma^f(f_\gamma, W_\gamma, \tau_2)\cap \Ou_{\gamma^\prime}^b(f_{\gamma^\prime},W_{\gamma^\prime}, -\tau_1) \cap
V^{-1}_\gamma(c) \cap
V^{-1}_{\gamma^\prime}(c)$.

The proof of 3. is similar to that of 1. and we omit it.

We now demonstrate 2. Suppose that $(x,\sigma)$ verifies
Assumption \ref{W2};  since $(x,\sigma)$
is precompact relative to $\Ou$, there is a compact set $B\subset
\Ou$ such that $x(t)\in B$ for all $t\ge 0$ and therefore
$(x(t),\sigma(t))$ belongs to the compact set $\dom(f)\cap (B
\times \Gamma) \subset \dom(W)$ for all $t\ge 0$. That 
$y(t)=W(x(t),\sigma(t))$ is weakly meagre follows from the fact that  $\int_0^\infty y(t)dt$ is
finite. Let us prove this last fact.

As $v(t)=V(x(t),\sigma(t))$ in nonincreasing and differentiable on each finite
interval $[a,b)$ in which  $\sigma$ is constant, (say $\sigma(t)=\gamma$ for all $t \in
[a,b)$), from \cite[Theorem 7.21]{Wheeden} we have that \[v(a)-v(b^-)\ge \int_{a}^b
(-\dot{v})(s)ds=-\int_a^b \nabla V_\gamma(x(s))f_\gamma(x(s))\;ds.\] Taking into
account (\ref{des1}) and that $v$ is nonincreasing on $[0,\infty)$, we then have
 that 
$$0\le \int_0^t y(s) ds \le V(x(0),\sigma(0))-V(x(t),\sigma(t))\;\; \forall t \geq 0.$$

As $(x(t),\sigma(t))$ evolves in the compact set $\dom(f)\cap (B\times \Gamma)$ for
all $t\ge 0$ and $V$ is bounded on compact subsets of $\dom(f)\cap(\Ou \times
\Gamma)$ then, for some $M\ge 0$, $|V(x(t),\sigma(t))|\le M$ for all $t\ge 0$. Thus $
\int_0^t y(s) ds \le 2M$ for all $t\ge 0$ and the integrability of $y$ follows.

Since $(x,\sigma)$ verifies the hypotheses of Theorem \ref{thHG} with $W$ instead of
$h$, $x\to \pi_1(M^*)$ where $M^*$ is the maximal weakly-invariant set w.r.t.
$\Ta[\tau_D,N_0]$ contained in $\dom(f)\cap W^{-1}(0)$.

In order to prove that
$\pi_1(M^*)\subset \bigcup_{\gamma, \gamma^\prime \in \Gamma} \left( \Ou_\gamma^f(f_\gamma, W_\gamma)\cap \Ou_{\gamma^\prime}^b(f_{\gamma^\prime},W_{\gamma^\prime})\right)$, we can
proceed as in the proof of 1.,  but using now the fact that for all $(\xi, \gamma) \in M^*$ there exists $(x^*, \sigma^*) \in \Ta[\tau_D, N_0]$ such that $0 = W_{\sigma^*(t)}(x^*(t)) = W_{\sigma^*(t^-)}(x^*(t))$ for all $t \in \R$ (being the last equality due to the continuity of $W$  on its domain). \qed
\begin{remark}  Theorem
\ref{convergence0} gives a more accurate result than Theorem 8 in \cite{Hespanha} in
the case when the switching signal $\sigma \in \Sa$ (instead of $\sigma \in
\So_{p-dwell}$ as is considered there). In fact, it can be shown that the hypotheses
of that theorem imply that the forward complete trajectory $(x,\sigma)$ of the linear
switched system $\dot{x}=A_{\sigma(t)}x$ is precompact and verifies Assumption
\ref{W2} with $V(\xi,\gamma)=\xi^TP_\gamma \xi$ and
$W(\xi,\gamma)=\xi^TC_{\gamma}^TC_\gamma\xi$. So, by applying Theorem
\ref{convergence0}, and taking into account 2. of Remark \ref{obs}, it results that
$x \to \cup_{\gamma \in \Gamma} \mathcal{U}_\gamma$, where $\mathcal{U}_\gamma$ is
the unobservable subspace of the pair $(C_\gamma,A_\gamma)$. On the other hand,
Theorem 8 in \cite{Hespanha} asserts that $x\to \mathcal{M}$, where $\mathcal{M}$ is
the smallest subspace which contains $\cup_{\gamma \in \Gamma} \mathcal{U}_\gamma$
and is $A_\gamma$-invariant for all $\gamma \in \Gamma$.
\end{remark}
\begin{remark} If in addition to the hypotheses of Theorem
\ref{convergence0}, we have that for some $x_e \in \cup_{\gamma
\in \Gamma} \Ou_\gamma$,  either for all $\gamma \in \Gamma$,  $\Ou^f_\gamma(f_\gamma,W_\gamma)\subset
\{x_e\}$ or for all $\gamma \in \Gamma$, $\Ou^b_\gamma(f_\gamma,W_\gamma)\subset
\{x_e\}$,
 then $x \to x_e$. The first conclusion of
Corollary 4.10 in \cite{Goebel} is a particular case of this result.

We note that, according to the particular geometry of each $\chi_\gamma$, it could happen that $\Ou^b_\gamma(f_\gamma,W_\gamma) \neq \Ou^f_\gamma(f_\gamma,W_\gamma)$ and even that 
one of those sets be void and the other one not.
\end{remark}

In what follows let for each $\gamma \in \Gamma$,  $E_\gamma=\{\xi \in
\chi_\gamma:\;f_\gamma(\xi)=0\}$   the set of equilibrium points of $f_\gamma$.

The following convergence result involves an ``ergodicity" condition on the
switching signals considered.
\begin{theorem} \label{ergodicconv} Suppose that $\Gamma$ is a finite set. Let
$(x,\sigma)$, with $\sigma \in \So_e \cap \Sd$, be a forward
complete trajectory of (\ref{ss}). Then the following holds:
\begin{enumerate}
 \item if $(x, \sigma)$ verifies Assumption \ref{W} and
 if for every $\gamma \in
\Gamma$, either $\Ou_{\gamma}^b(f_\gamma,W_\gamma)$ $=E_\gamma \cap \Ou_\gamma$ or $\Ou_{\gamma}^f(f_\gamma,W_\gamma)=E_\gamma \cap \Ou_\gamma$, then there exists $c\in \R$ such that $x \to \cap_{\gamma \in \Gamma}
(E_\gamma \cap V_\gamma^{-1}(c))$. If, in addition, for each $c\in \R$,
$\cap_{\gamma \in \Gamma} (E_\gamma \cap V_\gamma^{-1}(c))$ is a discrete set, then
$x \to x_e$ for some $x_e \in \cap_{\gamma \in \Gamma} (E_\gamma \cap \Ou_\gamma)$.
\item If $(x, \sigma)$ verifies Assumption \ref{Wf} and if for every $\gamma \in
\Gamma$, either $\Ou_{\gamma}^b(f_\gamma,W_\gamma)$ $=E_\gamma \cap \Ou_\gamma$ or 
$\Ou_{\gamma}^f(f_\gamma,W_\gamma)=E_\gamma \cap \Ou_\gamma$, then $x \to \cap_{\gamma \in \Gamma} (E_\gamma \cap \Ou_\gamma)$. If, in
addition, $\cap_{\gamma \in \Gamma} (E_{\gamma} \cap \Ou_\gamma)$ is a discrete set,
then $x \to x_e$ for some $x_e \in \cap_{\gamma \in \Gamma} (E_\gamma \cap
\Ou_\gamma)$
\end{enumerate}
\end{theorem}
{\bf Proof.} As $\sigma \in \So_e \cap \Sd$, then there exist $T>0$ and $\tau_D>0$ such that
$\sigma \in \So_e[T] \cap \Sd[\tau_D]$.

Let us  prove 1. first. Since $\So_e[T] \cap \Sd[\tau_D]$ has property {\bf P} and $(x,\sigma)$ verifies the hypotheses of Theorem \ref{thVG}, there exists $c\in\R$
such that $x \to \pi_1(M(c))$, where $M(c)$ is the maximal weakly-invariant set
w.r.t. $\T_e[T]\cap\Td[\tau_D]$ contained in $V^{-1}(c)\cap Z_V$. So, it suffices to
show that $\pi_1(M(c))\subset \cap_{\gamma \in \Gamma} (E_\gamma \cap
V_\gamma^{-1}(c))$.

Let $(\xi,\gamma)\in M(c)$;
then there exists a trajectory $(x^*,\sigma^*) 
\in \T_e[T]\cap \Td[\tau_D]$ such that $(x^*(0),\sigma^*(0))=(\xi,\gamma)$ and 
$(x^*(t),\sigma^*(t)) \in M(c)$ for all $t \in \R$. Then, reasoning as in the proof of Theorem \ref{convergence0},
$V_{\sigma^*(t^-)}(x^*(t)) \equiv V_{\sigma^*(t)}(x^*(t)) \equiv c$ 
and, from (\ref{des1}),
$W_{\sigma^*(t^-)}(x^*(t)) \equiv W_{\sigma^*(t)}(x^*(t))\equiv 0$. 

Taking into account that $x^*$
is continuous and that either 
$\mathcal{O}_{\gamma}^b(f_\gamma,W_\gamma)=E_\gamma \cap \Ou_\gamma$ or $\mathcal{O}_{\gamma}^f(f_\gamma,W_\gamma)=E_\gamma \cap \Ou_\gamma$, it follows that $x^*(t)\equiv \xi$ and that $\xi \in
E_{\sigma^*(t)}\cap \Ou_{\sigma^*(t)} \cap V_{\sigma^*(t)}^{-1}(c)=E_{\sigma^*(t)}\cap
V_{\sigma^*(t)}^{-1}(c)$ for all $t \in \R$. 

As $\sigma^*\in \So_e[T]$, then
$\sigma^*(\R)=\Gamma$ and, in consequence,  $\xi \in E_\gamma \cap V_\gamma^{-1}(c)$
for all $\gamma \in \Gamma$.

In the case that for every $c \in \R, \,\cap_{\gamma \in \Gamma} \left(E_{\gamma} \cap V_\gamma^{-1}(c)\right)$ is a discrete set, that $x\to x_e$ with $x_e \in \cap_{\gamma \in \Gamma}
(E_{\gamma} \cap \Ou_\gamma)$, follows from the facts that $x \to \Omega(x)$, that  $\Omega(x)$ is a nonempty
connected set and
that $\Omega(x) \subset \cap_{\gamma \in \Gamma} \left(E_{\gamma} \cap V^{-1}_\gamma (c) \right)$ for some $c \in \R$.

The proof of 2. is similar to that of 1., so we omit it. \qed

In the sequel we give sufficient conditions for the convergence to a given
equilibrium point $x_e$ of (\ref{ss}), i.e. a point $x_e$ that verifies $f_\gamma(x_e) = 0$ for all $\gamma \in
\Gamma$ such that $x_e \in \chi_\gamma$. We assume, without loss of generality, that
$x_e$ is the origin.

\begin{as} \label{f0}
$0$ is an equilibrium point of  (\ref{ss}). \end{as}
\begin{as} \label{uc} For every $\gamma \in \Gamma$ such
that $0 \in \chi_\gamma$, the initial value problem $\dot{x}=f_\gamma(x)$, $x(0)=0$
has a unique solution.
\end{as}

\begin{theorem} \label{convergence1} Suppose that assumptions \ref{f0}
and \ref{uc} hold and let $(x,\sigma)$ be a forward complete trajectory of
(\ref{ss}) with $\sigma \in \Sa$.
\begin{enumerate}
 \item If Assumption \ref{W2} is verified, $0 \in \Ou$ and the following holds
\begin{enumerate}
\item $\mathcal{O}_\gamma^f(f_\gamma,W_\gamma,\infty)\cap
\mathcal{O}_\gamma^b(f_\gamma,W_\gamma,\infty)\subset \{0\}$ for every $\gamma \in
\Gamma$, \item $\Ou_{\gamma}^b(f_{\gamma},W_{\gamma})\cap
\Ou_{\gamma^\prime}^f(f_{\gamma^\prime},W_{\gamma^\prime}) \subset \{0\},\quad
\forall \gamma\neq \gamma^\prime$,
\end{enumerate}
then $x \to 0$.

If $\Gamma$ is finite, the same holds if we suppose that $(x,\sigma)$ verifies
Assumption \ref{Wf} instead of Assumption \ref{W2}. \item If Assumption \ref{W} is
verified, $0 \in \Ou$, {\rm 1.(i)} holds and
\begin{enumerate}
\setcounter{enumi}{1}
 \item
$\Ou_{\gamma}^b(f_{\gamma},W_{\gamma})\cap
\Ou_{\gamma^\prime}^f(f_{\gamma^\prime},W_{\gamma^\prime})\cap
V_\gamma^{-1}(c) \cap V_{\gamma^\prime}^{-1}(c) \subset \{0\},\;\;
\forall \gamma\neq \gamma^\prime \in \Gamma, \;\; \forall c \in
\R$,
\end{enumerate}
then $x \to 0$.
\end{enumerate}
\end{theorem}
{\bf Proof.}
 Since $\sigma \in \Sa$, there exist $\tau_D>0$ and
$N_0 \in \N$ such that $\sigma \in \Sa[\tau_D,N_0]$.

We first prove 2. By using the same arguments as in the proof of the first part of
Theorem \ref{convergence0}, it follows that there exists $c \in \R$ such that $x \to
\pi_1(M(c))$, with $M(c)$ the maximal weakly-invariant set w.r.t. $\Ta[\tau_D,N_0]$
contained in $V^{-1}(c)\cap Z_V$. So, it suffices to show that $M(c) \subset \{0\}
\times \Gamma$.

Let $(\xi,\gamma)\in M(c)$; then there exists a trajectory $(x^*,\sigma^*) \in
\Ta[\tau_D,N_0]$ such that $(x^*(0),\sigma^*(0))=(\xi,\gamma)$ and
$(x^*(t),\sigma^*(t)) \in M(c)$ for all $t \in \R$. Once again, as in 
Theorem \ref{convergence0} we have that 
$V_{\sigma^*(t^-)}(x^*(t)) \equiv V_{\sigma^*(t)}(x^*(t)) \equiv c$ 
and $W_{\sigma^*(t^-)}(x^*(t)) \equiv W_{\sigma^*(t)}(x^*(t))\equiv 0$. 
  We will consider
two cases. \\
{\em Case 1.} $\sigma^*$ has no switching times, i.e. $\sigma^*(t)=\gamma$ for all
$t \in \R$. Then for every $t\in \R$, $x^*(t)\in \Ou_\gamma$ and
$W_{\gamma}(\varphi(t))=0$ and, since $x^*(t)\in \mathcal{O}_\gamma^f(f_\gamma,W_\gamma,\infty)\cap
\mathcal{O}_\gamma^b(f_\gamma,W_\gamma,\infty)\subset \{0\}$, $x^*(t)=0$ for all $t$. \\
{\em Case 2.} $\sigma^*$ has a switching time $t^*$. Then, there exists $\tau>0$
such that, if $\varphi(t)=x^*(t+t^*)$, $\hat{\gamma}=\sigma^*(t^{*-})$ and
$\gamma^\prime=\sigma^*(t^*)$, $\hat{\gamma} \neq \gamma^\prime$ and
\begin{enumerate}
\item $\varphi : [-\tau, 0] \rightarrow \Ou_{\hat{\gamma}}$ is solution of $\dot{z}=f_{\hat{\gamma}}(z)$ on $[-\tau,0]$ and $\varphi : [0, \tau] \rightarrow \Ou_{\gamma^\prime}$ is solution of
$\dot{z}=f_{\gamma^\prime}(z)$ on $[0,\tau]$ ; \item
$V_{\hat{\gamma}}(\varphi(t))=c$ on $[-\tau,0]$ 
and $V_{\gamma^\prime}(\varphi(t))=c$ on $[0,\tau]$;
\item  $W_{\hat{\gamma}}(\varphi(t))=0$
on $[-\tau,0]$ and $W_{\gamma^\prime}(\varphi(t))=0$ on $[0,\tau]$.
\end{enumerate}
Thus $x^*(t^*) = \varphi(0) \in \mathcal{O}_{\hat{\gamma}}^b(f_{\hat{\gamma}}, W_{\hat{\gamma}})\cap
\mathcal{O}_{\gamma^\prime}^f(f_{\gamma^\prime}, W_{\gamma^\prime})\cap
V_{\hat{\gamma}}^{-1}(c) \cap V_{\gamma^\prime}^{-1}(c)\subset \{0\}$.

That $x^*(0)=0$ follows from the fact that, due to assumptions \ref{f0} and
\ref{uc}, any initial value problem $\dot{z}=f_{\hat{\gamma}}(z)$, $z(0)=0$ has the
unique solution $z(t)\equiv 0$ when $0 \in \Ou_{\hat{\gamma}}$.

In order to prove 1. we note that in the case in which $\Gamma$ is finite and
$(x,\sigma)$ verifies Assumption \ref{Wf} then, due to Theorem \ref{theoremVG} (with
$\So^*=\Sa[\tau_D,N_0]$), it suffices to show that for any $\vec{c} \in \R^N$, the set
$M(\vec{c})$ of that theorem is a subset of $\{0\}\times \Gamma$. Since
the proof of this fact is similar to that of $M(c)\subset\{0\} \times \Gamma$ given
above, we omit it.

Suppose now that $(x,\sigma)$ verifies Assumption \ref{W2}. It follows from the
proof of Theorem \ref{convergence0} that $y(t)=W(x(t),\sigma(t))$ is weakly meagre. Since
$W$ is continuous on $\dom(f)\cap (\Ou \times \Gamma)$, $(x,\sigma)$ verifies the
hypotheses of Theorem \ref{thHG} and in consequence $x\to \pi_1(M^*),$ being $M^*$ the maximal weakly-invariant set w.r.t. $\Ta[\tau_D,N_0]$ contained in
$\dom(f)\cap W^{-1}(0)$. The proof of $M^* \subset \{0\} \times \Gamma$ in similar
to that of $M(c) \subset \{0\} \times \Gamma$, so we only delineate it. Let
$(\xi,\gamma)\in M^*$; then there exists $(x^*,\sigma^*)\in \Ta[\tau_D,N_0]$ such
that $(x^*(0),\sigma^*(0))=(\xi,\gamma)$ and, as in Theorem \ref{convergence0},  $W_{\sigma^*(t^-)}(x^*(t)) = W_{\sigma^*(t)}(x^*(t))=0$ for all
$t\in \R$. If $\sigma^*$ has no switching times, then $x^*(0)=0$ due to 1.(i).

If $\sigma^*$ has a switching time $t^*$ we considerer $\tau$, $\varphi$,
$\hat{\gamma}$ and $\gamma^\prime$ as in the proof of Case 2. above. Then $\varphi$
verifies 1. of that case, $W_{\hat{\gamma}}(\varphi(t))=0$ on $[-\tau,0]$
 and $W_{\gamma^\prime}(\varphi(t))=0$ on
$[0,\tau]$. Thus, from 1.(ii) we deduce that $x^*(t^*) = \varphi(0)=0$ and {\em a posteriori}
that $x^*(0)=0$.
\qed

When $\Gamma$ is finite and $\sigma$ belongs to $\Sd \cap \So^H$, hypothesis 2. of
Theorem \ref{convergence1} can be weakened as follows.

 Given a set-valued map $H:\Gamma \rightsquigarrow
\Gamma$, a finite sequence $\{\gamma_i\}_{i=1}^m \subset \Gamma$, $m\ge 3$, is a
{\em simple cycle of $H$} if $\gamma_1=\gamma_m$, $\gamma_{i+1}\in H(\gamma_i)$ for
all $i=1,\ldots,m-1$ and if $\gamma_i=\gamma_j$ and $i<j$ then $i=1$ and $j=m$.
\begin{theorem} \label{convergence2} Suppose that $\Gamma$ is finite,
that $H:\Gamma \rightsquigarrow \Gamma$ and that $(x,\sigma)$ is a forward complete
trajectory of (\ref{ss}) with $\sigma \in \Sd \cap \So^H$.
 Suppose in addition that assumptions \ref{f0}
and \ref{uc} hold.
\begin{enumerate}
\item If Assumption \ref{Wf} holds, $0 \in \Ou$ and
\begin{enumerate}
 \item
$\mathcal{O}_\gamma^f(f_\gamma,W_\gamma,\infty)\subset \{0\}$  for every
$\gamma \in \Gamma$ or $\mathcal{O}_\gamma^b(f_\gamma,W_\gamma,\infty)\subset \{0\}$ for every
$\gamma \in \Gamma$, \item for each simple cycle $\{\gamma_i\}_{i=1}^m$ of $H$ there
exists $j \in \{1,\ldots,m-1\}$ such that
\begin{eqnarray}\label{bf} \mathcal{O}_{\gamma_j}^b(f_{\gamma_j},W_{\gamma_j})\cap
\mathcal{O}_{\gamma_{j+1}}^f(f_{\gamma_{j+1}},W_{\gamma_{j+1}}) \subset
\{0\},\end{eqnarray}
\end{enumerate}
then $x \to 0$. \item The same conclusion as in {\rm 1.} holds if we replace Assumption
\ref{Wf} by Assumption \ref{W} and condition {\rm 1. (ii)}  by the weaker one:
\begin{enumerate}
\setcounter{enumi}{2}
\item  for every $c\in \R$ and for each simple cycle $\{\gamma_i\}_{i=1}^m$ of $H$
there exists $j \in \{1,\ldots,m-1\}$ such that
\begin{eqnarray}\label{bfw} \mathcal{O}_{\gamma_j}^b(f_{\gamma_j},W_{\gamma_j})\cap
\mathcal{O}_{\gamma_{j+1}}^f(f_{\gamma_{j+1}},W_{\gamma_{j+1}}) \cap
V_{\gamma_j}^{-1}(c)\cap V_{\gamma_{j+1}}^{-1}(c) \subset \{0\}.\end{eqnarray}
\end{enumerate}
\end{enumerate}
\end{theorem}
{\bf Proof.} As $\sigma \in \Sd \cap \So^H$, there exists $\tau_D>0$ such that $\sigma \in
\Sd[\tau_D]$.

Suppose that $(x,\sigma)$ verifies Assumption \ref{Wf}. Since $\Sd[\tau_D] \cap
\So^H$ has property {\bf P} and $(x,\sigma)$ verifies the hypotheses of Theorem
\ref{theoremVG}, there exists $\vec{c}\in \R^N$ such that $x \to \pi_1(M(\vec{c}))$, with $M(\vec{c})$ as
in that theorem (with $\T^* = \Td[\tau_D] \cap \T^H$). So, it suffices to show that $M(\vec{c})\subset \{0\} \times \Gamma$.

Let $(\xi,\gamma)\in M(\vec{c})$; then there exists a trajectory $(x^*,\sigma^*) 
\in \T_d[\tau_D]\cap \T^H$ such that $(x^*(0),\sigma^*(0))=(\xi,\gamma)$ and 
$(x^*(t),\sigma^*(t)) \in M(\vec{c})$ for all $t \in \R$. So, reasoning as in Theorem \ref{convergence0}, 
$V_{\sigma^*(t^-)}(x^*(t))=c_{\sigma^*(t^-)}$, $V_{\sigma^*(t)}(x^*(t))=c_{\sigma^*(t)}$  and $W_{\sigma^*(t^-)}(x^*(t)) = W_{\sigma^*(t)}(x^*(t))=0$ for all
$t\in \R$. We
distinguish
two cases:\\
{\em Case 1.} $\sigma^*$ has a finite number of switching times, $t_0 < t_1 <\cdots
< t_l$. Suppose first that for every $\gamma \in \Gamma$,  $\mathcal{O}_\gamma^f(f_\gamma,W_\gamma,\infty)\subset
\{0\}$  and let $\varphi(t)=x^*(t+t_l)$ for $t \ge 0$
and $\gamma_l=\sigma^*(t_l)$. Then $\varphi$ is a solution of
$\dot{z}=f_{\gamma_l}(z)$, $\varphi(t)\in \Ou_{\gamma_l}$ and 
$W_{\gamma_l}(\varphi(t))=0$ for all $t\ge 0$. In consequence $x^*(t_l)=0$, since
$x^*(t_l)=\varphi(0)\in \Ou_{\gamma_l}^f(f_\gamma,W_\gamma,\infty)\subset \{0\}$.

That $\xi=x^*(0)=0$, follows from the fact that for every $\gamma \in \Gamma$ such
that $0\in \Ou_\gamma$, the unique solution of the initial value problem
$\dot{z}=f_\gamma(z)$, $z(0)=0$ is
$z(t)\equiv 0$ .\\
In the case when for
every $\gamma \in \Gamma$, $\mathcal{O}_\gamma^b(f_\gamma,W_\gamma,\infty)\subset \{0\}$,  we proceed in a similar way, but
considering instead $\varphi(t)=x^*(t+t_0)$ for $t\le 0$.\\
{\em Case 2.} $\sigma^*$ has an infinite number of switching times. Since $\sigma^*
\in \So^H$, there exists a finite sequence of consecutive switching times
$\{t_{k}\}_{k=1}^m$, such that the sequence $\{\gamma_k\}_{k=1}^m$, with
$\gamma_{k}=\sigma^*(t_{k})$, is a simple cycle of $H$. By hypothesis there exists
an index $j \in \{1,\ldots,m-1\}$ for which (\ref{bf}) holds. For such $j$ we
consider the function $\varphi:[-\tau_D,\tau_D]\to \R^n$ defined by
$\varphi(t)=x^*(t+t_{j+1})$. Since $(x^*,\sigma^*) \in \Td[\tau_D]$ we have, for any
$0<\tau<\tau_D$, that
\begin{enumerate}
\item $\varphi : [-\tau, 0] \rightarrow \Ou_{\gamma_j}$ is solution of $\dot{z}=f_{\gamma_j}(z)$ on $[-\tau,0]$ and $\varphi : [0, \tau] \rightarrow \Ou_{\gamma_{j + 1}}$ is solution of
$\dot{z}=f_{\gamma_{j + 1}}(z)$ on $[0,\tau]$ ; 
\item 
 $W_{\gamma_j}(\varphi(t))=0$ on
$[-\tau,0]$ and $W_{\gamma_{j+1}}(\varphi(t))=0$ on $[0,\tau]$.
\end{enumerate}
Therefore $$\varphi(0) \in \mathcal{O}_{\gamma_j}^b(f_{\gamma_j}, W_{\gamma_j})\cap
\mathcal{O}_{\gamma_{j+1}}^f(f_{\gamma_{j+1}}, W_{\gamma_{j+1}})$$ and, by
(\ref{bf}), $x^*(t_{j+1})=\varphi(0)=0$. By using arguments similar to those of the
proof of case 1, we conclude that $\xi=x^*(0)=0$.

Since the proof of 2. is similar to that of 1., we omit it. \qed
\begin{remark} \label{obsV} 
It can be seen that Theorem \ref{convergence2} and Theorem \ref{convergence1} (supposing in Part 1. that $\Gamma$ is finite and that Assumption 5 holds) 
 remain valid if, instead of
Assumption \ref{uc}, we suppose that the function $V$ in assumptions \ref{W} and
\ref{Wf} verifies the following: for each $\gamma \in \Gamma$ such that $0 \in
\chi_\gamma$, 
$V^{-1}_\gamma(0) \cap \chi_\gamma = \{0\}.$
This condition is fulfilled when, for example,  $V_\gamma(\cdot)$ is positive
definite on $\chi_\gamma$ for every $\gamma$ such that $0 \in \chi_\gamma$.
\end{remark}

\subsection{Stability criteria}
Combining the convergence results already presented with well known sufficient
Lyapunov conditions for the local (global) stability of a family $\T$ of forward
complete trajectories of (\ref{ss}), we can derive some new local (global)
asymptotic stability criteria.

We recall that a family $\T$ of forward complete trajectories of (\ref{ss}) is
\begin{enumerate}
\item {\em locally uniformly stable} (LUS) if there exist a positive number $r>0$
and a function $\alpha:[0,r]\to \R$ of class $\kk$\footnote{
As usual, by a $\kk$-function we mean a function
$\alpha: [0, r] \rightarrow \R_{\ge 0}$ that is strictly increasing and 
continuous, and satisfies $\alpha(0)=0$. A $\ki$-function is one of class $\kk$ for which $r = + \infty$ and that is in addition unbounded.} such that for all $(x,\sigma)\in
\T$
$$ |x(t_0)|\le r \Rightarrow |x(t)|\le \alpha(|x(t_0)|) \quad \forall
t\ge t_0, \forall t_0\ge 0;$$ \item {\em globally uniformly stable} (GUS) if there
exists a function $\alpha:[0,\infty) \to \R$ of class $\kk_{\infty}$ such that for
all $(x,\sigma)\in \T$
$$ |x(t)|\le \alpha(|x(t_0)|) \quad \forall
t\ge t_0, \forall t_0\ge 0;$$ \item {\em locally asyptotically stable} (LAS) if it
is LUS and there exists $\eta>0$ such that for all $(x,\sigma) \in \T$ with
$|x(0)|<\eta$, $x\to 0$; \item {\em globally asyptotically stable} (GAS) is it is
GUS and for all $(x,\sigma) \in \T, x\to 0$.
\end{enumerate}

The different stability results that we present next,  require the introduction of
the following pair of functions.
\begin{definition}\label{wlf}
We say that a pair $(V,W)$ is a {\em weak Lyapunov pair} for the family $\T$ of
forward complete trajectories of (\ref{ss}) if
\begin{enumerate}
\item $V \in \V$, $0 \in \Ou$ and there exist functions $\alpha_1$ and $\alpha_2$ of
class $\kk$ such that
\begin{eqnarray}\label{k}
\alpha_1(|\xi|)\le V(\xi,\gamma)\le \alpha_2(|\xi|)\quad \forall \xi \in \Ou_\gamma,
\forall \gamma \in \Gamma.
\end{eqnarray}
\item $W:\dom(f) \cap (\Ou \times \Gamma) \to \R$ is such that  (\ref{des1}) holds with
$W_\gamma(\cdot)=W(\cdot,\gamma)$.
\item For every $(x,\sigma)\in \T$, the following is verified: \\
$x(t) \in \Ou\;\; \forall t \in [a,b] \subset [0, + \infty)\; \Rightarrow$
$v(t)=V(x(t),\sigma(t))$ is nonincreasing on $[a,b]$.
\end{enumerate}
We say that a pair $(V,W)$ is a {\em F-weak Lyapunov pair} for the family $\T$ of
forward complete trajectories of (\ref{ss}) if $V$ and $W$ satisfy {\rm 1.} and {\rm 2.} and the
following condition, which is weaker than {\rm 3.}
\begin{enumerate}
\setcounter{enumi}{3}
\item For every $(x,\sigma)\in \T$, the following holds: \\
$x(t) \in \Ou$ for all $t \in [a,b] \subset [0, + \infty)$ $\Rightarrow$ for every
$\gamma \in \Gamma$ $v(t)=V(x(t),\gamma)$ is nonincreasing on $[a,b]\cap
\sigma^{-1}(\gamma)$.
\end{enumerate}
\end{definition}
By using standard techniques (like those in \cite{Branicky}
or in \cite{Liberzon-book}) it is not hard to prove that the existence of a weak Lyapunov pair (or a F-weak
Lyapunov pair when $\Gamma$ is finite) for a family of trajectories $\T$ of
(\ref{ss}), implies that $\T$ is LUS and that it is GUS if, in addition, $\Ou=\R^n$
and $V$ is radially unbounded, i.e. there exist functions $\alpha_1$ and $\alpha_2$
of class $\kk_\infty$ such that (\ref{k}) holds.

\begin{theorem} \label{guas1} Suppose that Assumption \ref{f0} holds and let $\T$ be a family of forward complete trajectories of
(\ref{ss}) such that for every $(x,\sigma) \in \T$, $\sigma \in \Sa$. Then $\T$ is
LAS if one of the following conditions holds:
\begin{enumerate}
\item there exists a weak Lyapunov pair $(V,W)$ for $\T$ such that  the
restriction of $V$ to $\dom(f)\cap(\Ou \times \Gamma)$ is continuous and {\rm 1.(i)} and
{\rm 2.(ii)} of Theorem \ref{convergence1} hold. \item Assumption \ref{uc} holds and
there exists a weak Lyapunov pair $(V,W)$ for $\T$ such that  $W$ is
continuous and {\rm 1.(i)} and {\rm 1.(ii)} of Theorem \ref{convergence1} hold. \item $\Gamma$
is finite and there exists a F-weak Lyapunov pair $(V,W)$ for $\T$ such that {\rm 1.(i)} and {\rm 1.(ii)} of Theorem \ref{convergence1} hold.
\end{enumerate}
If, in addition, $\Ou=\R^n$ and $V$ is radially unbounded, then $\T$ is GAS.
\end{theorem}
\begin{theorem}\label{guas2} Suppose that $\Gamma$ is finite and
that Assumption \ref{f0} holds. Let $\T$ be a family of forward complete
trajectories of (\ref{ss}) such that for every $(x,\sigma) \in \T$, $\sigma \in
\Sd\cap \So^H$, with $H:\Gamma \rightsquigarrow \Gamma$. Then $\T$ is LAS if one of
the following holds.
\begin{enumerate}
 \item There exists a weak Lyapunov pair $(V,W)$ such that  {\rm 1.(i)} and
 {\rm 2.(iii)} of Theorem \ref{convergence2} hold.
 \item There exists a F-weak Lyapunov pair $(V,W)$ such that {\rm 1.(i)} and
 {\rm 1.(ii)} of Theorem \ref{convergence2} hold.
\end{enumerate}
If, in addition, $\Ou=\R^n$ and $V$ is radially unbounded, then $\T$ is GAS.
\end{theorem}
\begin{theorem}\label{guas2bis} Suppose that $\Gamma$ is finite and let $\T$ be a family of forward complete trajectories of
(\ref{ss}) such that for every $(x,\sigma) \in \T$, $\sigma \in \So_e\cap \Sd$.
Suppose that there exists a F-weak Lyapunov pair $(V,W)$ for $\T$ such that
for all
$\gamma \in \Gamma$, either $\mathcal{O}^b_{\gamma}(f_\gamma,W_\gamma)=E_\gamma \cap \mathcal{O}_\gamma$ or $\mathcal{O}^f_{\gamma}(f_\gamma,W_\gamma)=E_\gamma \cap \mathcal{O}_\gamma$  and that $\cap_{\gamma \in \Gamma} \left(E_\gamma \cap \Ou_\gamma \right)=\{0\}$.
 Then $\T$ is LAS. 

If, in addition, $\Ou=\R^n$
and $V$ is radially unbounded, then $\T$ is GAS.
\end{theorem}

{\em Proof of  Theorems \ref{guas1}, \ref{guas2} and \ref{guas2bis}}. Since the
hypotheses of the three theorems imply that $\T$ is locally uniformly stable (LUS),
we only need to prove that there exists $\eta >0$ such that for every $(x,\sigma)
\in \T$, $|x(0)|< \eta$ implies that $x\to 0$.

Since $\T$ is LUS, there exist $\eta>0$ and $\rho>0$ such that, for every
$(x,\sigma) \in \T$ with $|x(0)|<\eta$, $x(t)\in B=\{\xi \in \R^n:\:|\xi|\le
\rho\}\subset \Ou$ for all $t \ge 0$. Therefore $(x,\sigma) \in \T$ is precompact
relative to $\Ou$ whenever $|x(0)|<\eta$. Then, due to 1. of Remark \ref{obsV}, to
Theorem \ref{convergence1} in the case of Theorem \ref{guas1} and to Theorem
\ref{convergence2} in the case of Theorem \ref{guas2}, we have that $x\to 0$ for any
$(x,\sigma) \in \T$ such that $|x(0)| < \eta$.

In the case of Theorem \ref{guas2bis}, due to Theorem \ref{ergodicconv} we have that for every $(x,\sigma) \in \T$ such that $|x(0)| < \eta, \,
x \to \cap_{\gamma \in \Gamma}\left(E_\gamma \cap \Ou_\gamma\right)=\{0\}$. In consequence the local asymptotic stability of $\T$
follows.

When $\Ou=\R^n$ and $V$ is radially unbounded, we have that $\T$ is GUS. That $x\to
0$ for every $(x,\sigma) \in \T$ follows by using the fact that any trajectory of
$\T$ is precompact and the same arguments as above.
 \qed

\begin{remark} Theorem \ref{guas2bis} strengthens  Theorem 15 in \cite{Wang}
(which is the extension of the main result of
\cite{Cheng} to nonlinear switched systems ). In fact, the hypotheses of Theorem \ref{guas2bis} are weaker than
those of that theorem since, on one hand, even when restricted to the case
$V(\xi,\gamma)=V(\xi)$ and $\Ou=\R^n$ (as that theorem considers) the condition
$\cap_{\gamma \in \Gamma}E_\gamma=\{0\}$ is weaker than the hypothesis that $V$ is a
common joint Lyapunov function as is assumed in that work and, on the other hand, the
condition $\Ou_\gamma^f(f_\gamma,W_\gamma)=E_\gamma$ is weaker than the condition $M
\cap Z_\gamma=E_\gamma$ (with $Z_\gamma = \{\xi\,:\, W_\gamma(\xi) = 0\}$)
considered in \cite{Wang}, since it always holds that
$\Ou_\gamma(f_\gamma,W_\gamma) \subset M \cap Z_\gamma$ and sometimes the
inclusion is strict.
\end{remark}

From Theorem \ref{guas2bis} and Remark \ref{obs}.2. we can easily derive the
following result, that contains as a particular case Theorem 1 of \cite{Cheng}.
\begin{corollary}\label{final} Assume that $\Gamma$ is finite and that
$f_\gamma(\xi)=A_\gamma \xi$ with $A_\gamma \in \R^{n \times n}$ for all $\xi \in
\R^n$. Let $\T$ be a family of forward complete trajectories of (\ref{ss}) such that
for all $(x,\sigma) \in \T$, $\sigma \in \Sd\cap \So_e$. Suppose that there exists a
family of positive definite matrixes $\{P_\gamma, \gamma \in \Gamma\} \subset
\R^{n\times n}$ and a family of matrixes $\{C_\gamma, \gamma \in \Gamma\}$ such that
\begin{enumerate}
\item $P_\gamma A_\gamma+A_\gamma^T P_\gamma \le -C_\gamma^TC_\gamma$ for all
$\gamma \in \Gamma$; \item $v(t)=x^T(t)P_{\sigma(t)}x(t)$ is nonincreasing on
$[0,\infty)$ for all $(x,\sigma)\in \T$; \item for every $\gamma \in \Gamma$,
$\U_\gamma$, the unobservable subspace of the pair $(C_\gamma,A_\gamma)$, coincides
with $\ker(A_\gamma)$; \item $\cap_{\gamma \in \Gamma} \ker(A_\gamma)=\{0\}$.
\end{enumerate}
The, $\T$ is GAS.
\end{corollary}

\section{Conclusions}
In this paper we have obtained some invariance results for switched systems which satisfy
a dwell-time condition. These results enable us to study, in an unified way,
properties of bounded trajectories of switched systems whose switchings are
subjected not only to state-dependent constraints, but also to restrictions on the
accessibility from  each subsystem to other ones.

We also derived from these results some convergence and stability criteria. These
criteria involve observability-like conditions on functions which bound the
derivatives of some continuous functions that are nonincreasing along complete
trajectories of the switched systems.
 \begin{center}
{\bf ACKNOWLEDGMENTS}
\end{center}
 This work was partially supported by FONCyT Project 31255.

\appendix
\section{Proof of Lemma \ref{unif1}} 
The following lemma is used in the proofs of Lemma \ref{unif1} and Theorem
\ref{ginvariance}.
\begin{lem} \label{cp} Let $\{\sigma_k\}$ be a sequence of switching signals in
$\Sa[\tau_D,N_0]$ with $\tau_D>0$ and $N_0 \in \N$. Then there exist a subsequence
$\{\sigma_{k_l}\}$ and a switching signal $\sigma^* \in \Sa[\tau_D,N_0]$ such that
\begin{enumerate}
\item $\lim_{l \rightarrow \infty}\sigma_{k_l}(t)=\sigma^*(t)$ for
almost all $t\in \R$;\item for each $t \in \R$ there exists a sequence of times 
$\{r_l(t)\}$ such that $$\lim_{l\rightarrow \infty}r_l(t)=t,\;
 \lim_{l\rightarrow \infty}\sigma_{k_l}(r_l(t))=\sigma^*(t)\;\mbox{and}\;
\lim_{l\rightarrow \infty} \tau^1_{\sigma_{k_l}}(r_l(t))-r_l>0.$$
\end{enumerate}
\end{lem}
{\bf Proof.} Let $m \in \N$. By applying
Lemma A.1 in \cite{Mancilla-Garcia-scl06} to the sequence $\{\sigma_k(\cdot
-m)\}$,the thesis  holds with $[-m, + \infty)$ instead of $\R$. In addition, for any
$t\ge -m$ the sequence $\{r_l(t)\}$ can be chosen to satisfy the condition $r_l(t)>
-m$ for all $l$.

Due to this fact, for each  $m \in \N$ we can construct recursively a sequence of
positive integers $\{k^m_j\}_{j\in \N}$ such that:
\begin{itemize}
\item[a)] $\{k^{m+1}_j\}_{j\in \N}$ is a subsequence of $\{k^{m}_j\}_{j\in \N}$ for
all $m \in \N$; \item[b)] For each $m
\in \N$ there exists $\sigma^*_m \in \Sa[\tau_D,N_0]$ such that\\
1. $\lim_{j \rightarrow \infty}\sigma_{k^m_j}(t)=\sigma^*_m(t)$
for almost all $t\ge -m$;\\
2. for each $t \in [-m,\infty)$ there exists a sequence of times 
$\{r^m_j(t)\}_{j\in \N}$, with $r^m_j(t)> -m$ for all $j\in \N$, such that
$$\lim_{j\rightarrow \infty}r^m_j(t)=t,\; \lim_{j\rightarrow
\infty}\sigma_{k^m_j}(r^m_j(t))=\sigma^*_m(t)\;\mbox{and}\; \lim_{l\rightarrow
\infty} \tau^1_{\sigma_{k^m_j}}(r^m_j(t))-r^m_j(t)>0.
$$
In addition, if $m\ge 2$ and $t \ge -m-1$, $\{r^m_j(t)\}_{j\in \N}$ is a subsequence
of $\{r^{m-1}_j(t)\}_{j\in \N}$.
\end{itemize}

Let us  define $\{\sigma_{k_l}\}$ by
$\sigma_{k_l}=\sigma_{k^l_l}$ for all $l \in \N$, $\sigma^*$ by
$\sigma^*(t)=\sigma^*_l(t)$ if $t \in [-l,-l+1)$ with $l\ge 2$, and by
$\sigma^*(t)=\sigma_1(t)$ if $t \in [-1,\infty)$  and $\{r_l(t)\}$ by
$r_l(t)=r^l_l(t)$ for all $t \in \R$ and all $l \in \N$. Taking into account
that for all $t \ge -m$, $\{k^l_l\}_{l\ge m}$  and $\{r^l_l(t)\}_{l\ge m}$ are subsequences of $\{k^m_l\}_{l\in
\N}$ and of $\{r^m_l(t)\}_{l\in \N}$ respectively, the theorem follows. 
\qed

{\em Proof of Lemma \ref{unif1}}. Since it is clear that the three families of
switching signals verify 1. and 2. of definition \ref{P}, we only have to prove that
they also verify property 3. That this fact is true for $\Sa[\tau_D,N_0]$ follows
from Lemma \ref{cp}, and next we proceed to show that the other families of
switching signals also have this property.

Suppose that $\sigma_k \in \Sd[\tau_D] \cap \So^H$ for all $k$. Then, from Lemma
\ref{cp} and the fact that $\Sd[\tau_D]=\Sa[\tau_D,1]$ ,  there exist a subsequence
$\{\sigma_{k_l}\}$ and $\sigma^*\in \Sd[\tau_D]$ such that $\sigma_{k_l}\to
\sigma^*$ a.e. on $\R$. Let $t \in \R$ be a switching time of $\sigma^*$; since
$\sigma_{k_l}(\tau) \to \sigma^*(\tau)$ for almost all $\tau$, there exist $s < t <
s^\prime$ such that $\sigma_{k_l}(s) \to \sigma^*(s)$,
 $\sigma_{k_l}(s^\prime) \to
\sigma^*(s^\prime)$  and $s^\prime-s < \tau_D$. We note that
$\sigma^*(s^\prime)=\sigma^*(t)\neq \sigma^*(t^-)=\sigma^*(s)$.

Then, for $l$ large enough, $\sigma_{k_l}(s)\neq \sigma_{k_l}(s^\prime)$, and in
consequence, $\sigma_{k_l}$ has a unique switching time in the interval
$(s,s^\prime]$. Thus, $(\sigma_{k_l}(s),\sigma_{k_l}(s^\prime))\in \Gr(H)$ for $l$
large enough and, in consequence, since ${\rm
Graph}(H)$ is closed, 
$$(\sigma^*(t^-),\sigma^*(t))=(\sigma^*(s),\sigma^*(s^\prime))\in \Gr(H).$$ 
 Therefore $\sigma^*\in \Sd[\tau_D]\cap \So^H$.

Suppose now that $\{\sigma_k\} \subset \Sd[\tau_D] \cap \So_e[T]$. By the same
arguments of the previous case, there exist a subsequence $\{\sigma_{k_l}\}$ and
$\sigma^*\in \Sd[\tau_D]$ such that $\sigma_{k_l}\to \sigma^*$ a.e. on $\R$.

Let $t_0 \geq 0, \, \gamma \in \Gamma$ and $\epsilon > 0$ be fixed, and suppose that
\begin{eqnarray}
 \label{lemaerg1}
{\sigma^*}^{-1}(\gamma) \cap[t_0, t_0 + T+\epsilon] = \emptyset.
\end{eqnarray}
Let $\epsilon^\prime > 0$ such that $\epsilon^\prime < \min\{\epsilon/2, \tau_D/2\}$
and
\begin{eqnarray*}
 &&\lim_{k \rightarrow \infty} \sigma_k(t_0 + \epsilon^\prime) = \sigma^*(t_0 + \epsilon^\prime) = \gamma_0 \neq \gamma, \\
 &&\lim_{k \rightarrow \infty} \sigma_k(t_0 + T + \epsilon - \epsilon^\prime) = \sigma^*(t_0 + T + \epsilon - \epsilon^\prime) = \gamma_1 \neq \gamma.
\end{eqnarray*}
Then, since $\Gamma$ is a finite set, there exists $K_0 \in \N$ such that
\begin{eqnarray}
 \label{lemaerg2}
\sigma_k(t_0 + \epsilon^\prime) = \gamma_0\; \mbox{and}\; \sigma_k(t_0 + T +
\epsilon - \epsilon^\prime) = \gamma_1\quad \forall k\geq K_0.
\end{eqnarray}
Let $I = [t_0 + \epsilon^\prime, t_0 + T + \epsilon - \epsilon^\prime]$; as the
length of $I$ is greater than $T$ and (\ref{lemaerg2}) holds, then for each $k \geq
K_0$ there exists $t_k \in I$ such that $\sigma_k(t_k) = \gamma$ and
$\sigma_k(t_k^-) \neq \gamma$. From the compactness of $I$,  there exists a subsequence $\{t_{k_l}\} \subset
\{t_k\}$ that converges to, say, $t_\gamma \in I$.

Let $\epsilon^{\prime \prime} \in(0,\epsilon^\prime/2)$ and $L_0 \in \N$ such that
$t_{k_l} \in (t_\gamma -\epsilon^{\prime \prime}, t_\gamma + \epsilon^{\prime
\prime})$ for every $l \geq L_0$. Since for all $l \geq L_0$, $\sigma_{k_l } (s) =
\gamma \, \forall s \in [t_{k_l}, t_{k_l} + \tau_D)$, then for those $l$'s
$\sigma_{k_l } (s) = \gamma \, \forall s \in [t_\gamma + \epsilon^{\prime \prime},
t_\gamma + \epsilon^{\prime}) \subset[t_0, t_0 + T + \epsilon]$.
 Hence, there exists
 $t \in [t_\gamma + \epsilon^{\prime \prime}, t_\gamma + \epsilon^{\prime})$
 such that $\gamma = \lim_{l \rightarrow \infty} \sigma_{k_l}(t) =
 \lim_{k \rightarrow \infty} \sigma_k(t) =  \sigma^*(t)$,
 which contradicts (\ref{lemaerg1}). In consequence for every
 $\epsilon > 0$, every $t_0 \geq 0$ and every $\gamma \in \Gamma$,
 ${\sigma^*}^{-1}(\gamma) \cap[t_0, t_0 + T + \epsilon] \neq \emptyset$.

Suppose now that for certain $t_0 \geq 0$ and  $\gamma \in \Gamma$, 
${\sigma^*}^{-1}(\gamma) \cap[t_0, t_0 + T] = \emptyset$; then $\sigma^*(t_0 + T)
\neq \gamma$ and since $\sigma^*$ in right-continuous, there exists $\epsilon > 0$
such that $\sigma^*(t_0 + T + s) \neq \gamma \,\forall s \in [0, \epsilon]$ and then
${\sigma^*}^{-1}(\gamma) \cap[t_0, t_0 + T + \epsilon] = \emptyset$ which is a
contradiction. In consequence ${\sigma^*}^{-1}(\gamma) \cap[t_0, t_0 + T] \neq
\emptyset$ for every $t_0 \geq 0$ and every $\gamma \in \Gamma$, and the lemma
follows. \qed

\section{Proof of Theorem \ref{ginvariance}}
The next result, that shows that under suitable hypotheses certain families of
trajectories of system (\ref{ss}) enjoy a certain kind of sequential compactness, is
used in the proof of Theorem \ref{ginvariance}.

We say that a sequence $\{(x_k,\sigma_k)\}$ of trajectories of (\ref{ss}) is {\em uniformly 
precompact} if there exists a compact set $B\subset \R^n$ such that $x_k(t) \in B$
for all $t \in  I_{x_k}$.

Consider next the family $\mathcal{I}$ composed by the intervals $I=[a,b)$ or
$I=(-\infty,b)$ with $b$ finite or $b=+\infty$.

\begin{prop}  \label{SC} Let $\{(x_k,\sigma_k)\}$ be a uniformly precompact sequence of
trajectories of (\ref{ss}) such that for a certain interval $I \in \mathcal{I}$ and
certain $\tau_D>0$ and $N_0\in \N$, $\sigma_k \in \Sa[\tau_D,N_0]$ and $I_{x_k}=I$
for all $k$. Then, there exist a subsequence $\{(x_{k_l},\sigma_{k_l})\}$ and a
trajectory $(x^*,\sigma^*)$ of (\ref{ss}) with $I_{x^*}=I$ and $\sigma^*\in
\Sa[\tau_D,N_0]$ such that
\begin{enumerate}
\item $x_{k_l}\to x^*$ uniformly on compact subsets of $I$; \item $\sigma_{k_l} \to
\sigma^*$ a.e. on $\R$.
\end{enumerate}
\end{prop}
Proposition \ref{SC} makes use of the assertions of the following lemma.
\begin{lem} \label{cp2} Let $\{(x_k,\sigma_k)\}$ be a sequence
such that for all $k \in \N$, $x_k:[a,b]\to \R^n$ is absolutely continuous,
$\sigma_k:[a,b]\to \Gamma$ and $\dot{x}(t)=f(x_k(t),\sigma_k(t))$ a.e. on $[a,b]$.
Suppose that there exist a compact set $B \subset \R^n$ such that $x_k([a,b])
\subset B$ for all $k$,  and a function  $\sigma :[a,b]\to \Gamma$ such that
$\lim_{k\to \infty}\sigma_k(t)=\sigma(t)$ for almost all $t\in [a,b]$. Then there exist a
subsequence $\{x_{k_j}\}$ of $\{x_k\}$ and an absolutely continuous function
$x:[a,b]\to \R^n$ such that $x_{k_j} \to x$ uniformly on $[a,b]$ and
$\dot{x}(t)=f(x(t),\sigma(t))$ a.e. on $[a,b]$.
 \end{lem}
 {\bf Proof.} Since for all $k \in \N$, $(x_k(t),\sigma_k(t))\in \dom(f)\cap (B\times \Gamma), \dot{x}_k(t)=f(x_k(t),\sigma_k(t))$ a.e. $[a,b]$ and $f$ is
continuous on the compact set $\dom(f)\cap (B\times \Gamma)$, there exists $M>0$
such that for all $k \in \N$, $|\dot{x}_k(t)|\le M$ for almost all $t \in [a,b]$.
 In consequence $\{x_k\}$ is uniformly Lipschitz, since $|x_k(t)-x_k(s)|\le M |t-s|$ for all $t,s\in [a,b]$ and all
 $k$, and therefore equicontinuous.

Then, due to Arzel$\grave{{\rm a}}$-Ascoli Theorem, there
 exist a subsequence $\{x_{k_j}\}$ and a continuous function $x:[a,b]\to
 \R^n$ such that $x_{k_j}\to x$ uniformly on $[a,b]$.

We note that $(x(t),\sigma(t)) \in \dom(f)$ a.e. on $[a,b]$
 since for every $t \in [a,b]$ such that $\lim_{k \to \infty}\sigma_k(t)=\sigma(t)$ and  that $(x_k(t),\sigma_k(t))\in \dom(f)$ for all $k\in \N$, 
 $(x(t),\sigma(t)) \in \dom(f)$. 
 
Given that
 $\lim_{j \to \infty} \dot{x}_{k_j}(t)=\lim_{j \to \infty}
 f(x_{k_j}(t),\sigma_{k_j}(t))=f(x(t),\sigma(t))$ for almost all
 $t \in [a,b]$ and that $\{\dot{x}_{k_j}\}$ is majorized by the constant
 function $\rho(t)=M$ on $[a,b]$, by  applying Lebesgue's Dominated
 Convergence Theorem we have that
\begin{eqnarray*}
x(t)&=&\lim_{j\to \infty} \left ( x_{k_j}(a)+\int_a^t
 \dot{x}_{k_j}(s)ds\right )\\
 &=&x(a)+\int_a^t f(x(s),\sigma(s))\:ds,
 \end{eqnarray*}
 and the lemma follows.
 \qed

{\em Proof of Proposition \ref{SC}} Due to Lemma \ref{cp} we can suppose without
loss of generality that there exists 
  a switching signal $\sigma^* \in \Sa[\tau_D,N_0]$ such that
 $\lim_{k\to \infty}\sigma_k(t)=\sigma^*(t)$
 for almost all $t \in \R$.

 Consider a sequence of intervals $\{[a_k, b_k]\}$ such that for all $k \in \N$,
 $[a_k, b_k] \subset [a_{k + 1},
b_{k + 1}] \subset I$ and $\cup_{k \in \N} [a_k, b_k] = I$ and let a compact set $B$
such that $x_k(I)\subset B$ for all $k$. Since $x_k([a_1, b_1]) \subset B$, it
follows from Lemma \ref{cp2} that there exist a subsequence $\{(x_k^1,\sigma_k^1)\}$
and an absolutely continuous function $x^1 : [a_1, b_1] \rightarrow \R^n$ such that
$x_k^1 \rightarrow x^1$ uniformly on $[a_1, b_1]$ and
$\dot{x}^1(t)=f(x^1(t),\sigma^*(t))$ a.e. $[a_1,b_1]$. Applying Lemma \ref{cp2} to
the sequence $\{(x_k^1,\sigma_k^1)\}$ on the interval $[a_2,b_2]$ we can assure the
existence of a subsequence $\{(x_k^2,\sigma_k^2)\}$ of $\{(x_k^1,\sigma_k^1)\}$ and
an absolutely continuous function $x^2: [a_2,b_2]\to \R^n$ such that $x_k^2
\rightarrow x^2$ uniformly on $[a_2, b_2]$ and $\dot{x}^2(t)=f(x^2(t),\sigma^*(t))$
a.e. $[a_2,b_2]$. Proceeding in this way, it follows that for each $l \in \N$ we can
construct recursively a subsequence $\{(x_k^l,\sigma_k^l)\}$ of $\{(x_k,\sigma_k)\}$
such that:
\begin{enumerate}
\item $\{x^{l+1}_k\}$ is a subsequence of $\{x^{l}_k\}$ for all $l \in \N$; \item
for each $l \in \N$ there exists an absolutely continuous function $x^l : [a_l, b_l]
\rightarrow \R^n$ such that $x_k^l \rightarrow x^l$ uniformly on $[a_l, b_l]$ and
$\dot{x}^l(t)=f(x^l(t),\sigma^*(t))$ a.e. $[a_l,b_l]$.
\end{enumerate}
Let now the subsequence $\{(x^*_l,\sigma_l^*)\}$ of $\{(x_k,\sigma_k)\}$ defined for
each $l \in \N$ as $x^*_l(t) = x^l_l(t)$ for all $t \in I$ and
$\sigma^*_l(t)=\sigma_l^l(t)$ for all $t \in \R$ and let $x^* : I \rightarrow \R^n$
be the function given by $x^*(t) = x^l(t)$ for all $t \in [a_l, b_l]$ and all $l \in
\N$. Clearly, $x^*$ is locally absolutely continuous on $I$ and we show next that
$x^*_l \rightarrow x^*$ uniformly on compact subsets of $I$.

Let $K \subset I$ be compact; then there exists $j \in \N$ such that $K \subset
[a_j, b_j]$. Since for all $l \geq j, \; x^*_l \in \{x^j_k\}$, it follows that
$x^*_l \rightarrow x^j$ uniformly on $[a_j, b_j]$.  The uniform convergence of
$\{x^*_l\}$ to $x^*$ follows from the fact that $x^*(t) = x^j(t)$ for all $t \in
[a_j, b_j]$.

Let now for all $l \in \N$, $k_l \in \N$ such that
$(x_{k_l},\sigma_{k_l})=(x^*_l,\sigma_l^*)$. Then the subsequence
$\{(x_{k_l},\sigma_{k_l})\}$ verifies 1. and 2. of the thesis.

The proof finishes by showing that $(x^*,\sigma^*)$ is a trajectory of (\ref{ss}).
By construction we have that $\dot{x}^*(t)=f(x^*(t),\sigma^*(t))$ for almost all $t
\in I$. It remains to prove that $(x^*(t),\sigma^*(t))$ belongs to $\dom(f)$ for all
$t \in I$.

Let $t \in I$; then there exists a sequence of times $\{t_k\}$, such that: i) for
all $k$, $t_k \in I$, $t_k > t$, $(x^*(t_k),\sigma^*(t_k)) \in \dom(f)$ and
$\sigma^*(t_k)=\sigma^*(t)$, ii) $\lim_{k\to \infty}t_k=t$.

Hence, $(x^*(t_k),\sigma^*(t_k))\to (x^*(t),\sigma^*(t))$ and
$(x^*(t),\sigma^*(t))\in \dom(f)$ since $\dom(f)$ is closed.
 \qed

From Proposition \ref{SC} and by  arguments similar to those used in its proof, we
can derive the following
\begin{prop} \label{SC2} Let $\tau_D>0$ and $N_0\in \N$ and
let $\{(x_k,\sigma_k)\}$ be an uniformly precompact sequence of trajectories of
(\ref{ss}) such that $\sigma_k \in \Sa[\tau_D,N_0]$ and $[-k,\infty) \subset
I_{x_k}$ for all $k$. Then, there exist a subsequence $\{(x_{k_l},\sigma_{k_l})\}$
and a trajectory $(x^*,\sigma^*)\in \Ta[\tau_D,N_0]$ such that
\begin{enumerate}
\item $x_{k_l}\to x^*$ uniformly on compact subsets of $\R$; \item $\sigma_{k_l} \to
\sigma^*$ a.e. on $\R$.
\end{enumerate}
\end{prop}

{\em Proof of Theorem \ref{ginvariance}}. Without loss of generality we can suppose
that $I_x=[0, + \infty)$.
Let $(\xi,\gamma) \in \Omega^\sharp(x,\sigma)$; then there exists a 
strictly increasing and unbounded sequence $\{s_k\}$ which verifies 1. and 2. of
Definition \ref{omega}. Let $\sigma_k(\cdot)=\sigma(\cdot+s_k)$ and
$x_k(\cdot)=x(\cdot+s_k)$ and note that $\sigma_k \in \So^*$ for all $k$ since
$\So^*$ has property {\bf P}. As
 the sequence $\{(x_k,\sigma_k)\}$ is uniformly precompact, 
 $I_{x_k}=[-k,\infty)$ and, due to 1. of Definition \ref{P}, $\So^*\subset\Sa[\tau_D,N_0]$ for some $\tau_D>0$ and some $N_0\in\N$, then from   Proposition \ref{SC2} and taking into account 3. of Definition \ref{P}, there exist
 a subsequence $(x_{k_l},\sigma_{k_l})$ and a trajectory
$(x^*,\sigma^*) \in \T^*$ such that $\{x_{k_l}\}$ converges to $x^*$ uniformly on
compact subsets of $\R$ and $\{\sigma_{k_l}\}$ converges to $\sigma^*$ a.e. on $\R$.

Due to Lemma \ref{cp}, we can assume without loss of generality that
$\{\sigma_{k_l}\}$ verifies condition 2 of that lemma.

The proof ends provided that we show that $(x^*(0),\sigma^*(0))=(\xi,\gamma)$ and that 
$(x^*(t),\sigma^*(t))\in \Omega^\sharp(x,\sigma)$ for all $t \in \R$.

Let us prove first that $(x^*(0),\sigma^*(0))=(\xi,\gamma)$. From the fact that 
$x_k(0)=x(s_k)$, from 2. of Definition \ref{omega} and due to the convergence of $\{x_{k_l}\}$ to 
$x^*$, it holds that $x^*(0)=\xi$.

According to 2. of Lemma \ref{cp}, there exists a sequence $\{r_l\} \subset \R$ such 
that $\lim_{l \rightarrow \infty}r_l=0$, $\lim_{l \rightarrow 
\infty}\tau^1_{\sigma_{k_l}}(r_l)-r_l>0$ and $\lim_{l \rightarrow 
\infty}\sigma_{k_l}(r_l)=\sigma^*(0)$. From item 1. of Definition \ref{omega} and the 
fact that $\tau^1_{\sigma_k}(0)=\tau_\sigma^1(s_k) - s_k$ for all $k \in \N$, it follows 
that $\lim_{l \rightarrow \infty} \tau^1_{\sigma_{k_l}}(0)>0$. Then $r_l< 
\tau^1_{\sigma_{k_l}}(0)$ for $l$ large enough and therefore, $\sigma^*(0)=\lim_{l 
\rightarrow \infty}\sigma_{k_l}(r_l)=\lim_{l \rightarrow \infty}\sigma_{k_l}(0)=\lim_{l 
\rightarrow \infty}\sigma(s_{k_l})=\gamma$.

In order to prove that $(x^*(t),\sigma^*(t)) \in \Omega^\sharp(x,\sigma)$ for all $t
\in \R$, $t\neq 0$, let one such time $t$ and let $\{r_l(t)\}$ be a sequence as in
2. of Lemma \ref{cp}. In order to simplify the notation, we will write $r_l$ instead
of $r_l(t)$.

Consider the 
unbounded sequence $\{s^\prime_l\}$, defined by $s^\prime_l=r_l+s_{k_l}$, which we 
assume, without loss of generality, strictly increasing. Since
$\tau^1_{\sigma}(s^\prime_l)=\tau^1_{\sigma_{k_l}}(r_l) + s_{k_l}$, we have that 
$\lim_{l \rightarrow \infty} \tau^1_{\sigma}(s^\prime_l)-s^\prime_l>0$. So 
$\{s^\prime_l\}$ satisfies condition 1. of Definition \ref{omega}.

That $\lim_{l \rightarrow 
\infty}\sigma(s^\prime_l)=\lim_{l \rightarrow \infty}\sigma_{k_l}(r_l)=\sigma^*(t)$
follows from 2. of Lemma \ref{cp}. From the uniform convergence of $\{x_{k_l}\}$
to $x^*$ on compact sets and the continuity of $x^*$ we have that
$\lim_{l \rightarrow \infty} x(s^\prime_{l})=\lim_{l \rightarrow 
\infty}x_{k_l}(r_l)=x^*(t)$. Hence $(x(s^\prime_{l}),\sigma(s^\prime_{l})) \rightarrow 
(x^*(t),\sigma^*(t))$ as $l \rightarrow \infty$ and in consequence $(x^*(t),\sigma^*(t)) \in 
\Omega^\sharp(x,\sigma)$.\qed

\end{document}